\documentclass[12pt,oneside]{amsart}
\usepackage{amsmath,amssymb}
\makeatletter \renewcommand{\@biblabel}[1]{\hfill#1.}
\@addtoreset{equation}{section} \makeatother

\title{Hypersurface exceptional singularities}
\author{Shihoko ISHII}
\author{Yuri Prokhorov}
\address{Shihoko Ishii: Department of Mathematics,
Tokyo Institute of Technology, Oh-Okayama, Meguro, Tokyo, Japan}
\address{Yuri Prokhorov: Department of Mathematics,
Tokyo Institute of Technology, Oh-Okayama, Meguro, Tokyo, Japan \&
Moscow State University, Moscow, Russia }
\newcommand{\bZ}{{\mathbb Z}}
\newcommand{\bR}{{\mathbb R}}
\newcommand{\bN}{{\mathbb N}}
\newcommand{\bP}{{\mathbb P}}
\newcommand{\bC}{{\mathbb C}}
\newcommand{\bQ}{{\mathbb Q}}
\newcommand{\bW}{{\mathbb W}}
\newcommand{\bCn}{{{\mathbb C}^n}}
\newcommand{\bCnp}{{{\mathbb C}^n(\mathbf p)}}

\newcommand{\fp}{{f_{\mathbf p}}}
\newcommand{\fq}{{f_{\mathbf q}}}
\newcommand{\ba}{{\mathbf a}}
\newcommand{\bp}{{\mathbf p}}
\newcommand{\bq}{{\mathbf q}}
\newcommand{\bl}{{\mathbf 1}}
\newcommand{\be}{{\mathbf e}}
\newcommand{\Xp}{{X(\mathbf p)}}
\newcommand{\Xq}{{X(\mathbf q)}}
\newcommand{\D}{{\Delta}}
\newcommand{\Dp}{{D_{\bp}}}
\newcommand{\Dq}{{D_{\bq}}}
\newcommand{\tnd}{{T_{N}(\Delta)}}
\newcommand{\mult}{\operatorname{mult}}
\newcommand{\dif}{\operatorname{Diff}}

\newcommand{\up}[1]{\left\lceil #1\right\rceil}
\newcommand{\De}{{\varDelta}}
\newcommand{\Supp}{\operatorname{Supp}}
\newcommand{\OOO}{{\mathcal{O}}}
\newcommand{\qq}{\mathbin{\sim_{\scriptscriptstyle{\mathbb{Q}}}}}
\newcommand{\ep}{\varepsilon}
\newcommand{\var}{\varphi}
\renewcommand{\bar}[1]{\overline{#1}}
\newcommand{\Ga}{\varGamma}
\newcommand{\const}{c}
\newcommand{\compl}{\operatorname{compl}}
\newcommand{\rk}{\operatorname{rk}}
\newcommand{\Pic}{\operatorname{Pic}}
\newcommand{\Sing}{\operatorname{Sing}}
\newcommand{\lcm}{\operatorname{lcm}}
\newcommand{\down}[1]{\left\lfloor #1\right\rfloor}
\newcommand{\fr}[1]{\left\{ #1 \right\}}
\renewcommand{\emptyset}{\varnothing}
\newcommand{\eref}[1]{{\rm{(\ref{#1})}}}
\newcommand{\x}[4]{[#1,#2,#3,#4]}
\newcommand{\y}[4]{y_1^{#1}+y_2^{#2}+y_3^{#3}+y_4^{#4}}
\newtheorem{thm}[subsection]{Theorem}
\newtheorem{lem}[subsection]{Lemma}
\newtheorem{cor}[subsection]{Corollary}
\newtheorem{prop}[subsection]{Proposition}

\newtheorem*{proposition*}{Proposition}
\newtheorem*{conjecture*}{Conjecture}

\newtheorem{corollary1}[subsubsection]{Corollary}

\theoremstyle{definition}
\newtheorem{defn}[subsection]{Definition}
\newtheorem{say}[subsection]{}
\newtheorem{exmp}[subsection]{Example}

\newtheorem{rem}[subsection]{Remark}
\newtheorem*{definition*}{Definition}

\theoremstyle{remark}

\newtheorem*{remark*}{Remark}

\newcounter{no}
\newcommand{\no}{\addtocounter{no}{1}\arabic{no}}

\begin{document}
\begin{abstract}
This paper studies hypersurface exceptional singularities in
$\bCn$ defined by non-degenerate function.
For each canonical hypersurface singularity, there exists a
weighted homogeneous singularity such that the former is
exceptional if and only if the latter is exceptional.
So we study the weighted homogeneous case and prove that the
number of weights of weighted homogeneous exceptional
singularities are finite.
Then we determine all exceptional singularities of the Brieskorn
type of dimension $3$.
\end{abstract}
\maketitle
\section{Introduction}
The notion of exceptional singularities was introduced and its
$2$-dimensional examples were given by Shokurov in \cite[\S
5]{shokurov2}.
At the beginning, there was no non-trivial example of higher
dimensional exceptional singularities.
But in \cite{m-p} the first higher dimensional non-trivial
examples of exceptional singularities are found and then in
\cite{m-p2} $3$-dimensional exceptional quotient singularities are
characterized.
The advantage of distinguishing exceptional and non-exceptional
singularities is as follows:
\begin{enumerate}
\item[(1)]
for a non-exceptional singularity, the linear system $|-mK_{X}|$
is to have a ``good'' member for some small $m$ (actually, we can
take $m\in \{1,2\}$ in the $2$-dimensional case and $m\in
\{1,2,3,4,6\}$ in the $3$-dimensional case, see \cite[5.2.3 and
Theorem~5.6]{shokurov2}, \cite[Theorem 4.9]{prokhorov} and also
Proposition~\ref{prokhorov});
\item[(2)]
exceptional singularities are to be classified.
\end{enumerate}
In this paper, we try to determine hypersurface exceptional
singularities defined by non-degenerate function.
Since an exceptional singularity is log-canonical, under our
situation it must be either strictly log-canonical or canonical.
For the trivial case, that is, the strictly log-canonical case,
the singularity is exceptional if and only if it is purely
elliptic of type $(0, d-1)$, where $d$ is the dimension of the
singularity (\cite{ishii}).
(In \cite{ishii} these terminologies are defined only for an
isolated singularity, but these are naturally extended to the case
that the log-canonical singular locus is isolated).
If $d=2$, then it is a simple elliptic singularity and there are
3-types $\tilde E_{6}$, $\tilde E_{7}$ and $\tilde E_{8}$ (see
\cite{saito}). If $d=3$ then it is a simple $K3$-singularity (c.f.
\cite{ishii-watanabe}), and isolated weighted homogeneous ones are
classified into 95-types in \cite{yonemura}.
Here we study the other case in which the singularity is
canonical.
In the $2$-dimensional case, canonical singularities are
hypersurface and there is a plt blow-up which is given by a
weighted blow-up. Moreover, the class of exceptional canonical
$2$-dimensional singularities is bounded: $E_{6}$, $E_{7}$,
$E_{8}$.
We generalize these facts to the higher dimensional case.
A criterion for a canonical singularity to be an exceptional
singularity is obtained in \cite{prokhorov} by means of plt
blow-up.
In order to make use of it, we have to construct a plt-blow-up
first.
In \S 3, we prove that there exists a weighted blow-up which gives
a plt-blow-up. We also prove that for every canonical singularity
defined by a non-degenerate function there exists such a
singularity defined by a weighted homogeneous polynomial such that
the former is exceptional if and only if the latter is
exceptional.
So we can reduce the problem into the weighted homogeneous case.
In \S 4, we prove the finiteness of the set of all weights of
weighted homogeneous hypersurface exceptional singularities of
fixed dimension.
In \S 5, we determine all exceptional singularities of the
Brieskorn type of dimension $3$.

\section{Preliminaries}
All varieties are defined over $\bC$. We use terminologies lc,
klt, plt, dlt defined in \cite{kollar}, \cite{uta}.
We denote a germ of a singularity $x$ in $X$ by $(X,x)$.
Let $D$ be a normal subvariety of codimension one on a normal
variety $X$. Assume $D$ is Gorenstein in codimension one. Then by
\cite[16.5]{uta} there exists a $\bQ$-Weil divisor $\dif_{D/X}(0)$
on $D$ such that:
\[ K_{D}+\dif_{D/X}(0)=K_{X}+D|_{D}. \]
This divisor $\dif_{D/X}(0)$ is called the \textsl{different}. If
there is no possibility of confusion, this is written as
$\dif_{D}(0)$.

\begin{defn}
Let $(X,x)$ be a normal singularity and $D=\sum d_{i}D_{i}$ a
boundary on $X$ such that $K_{X}+D$ is log-canonical. The pair
$(X,D)$ is said to be \textsl{exceptional} if there exists at most
one exceptional divisor $E$ over $X$ with discrepancy $a(E,D)=-1$.
A log-canonical singularity $(X,x)$ is said to be
\textsl{exceptional} if every lc pair $(X,D)$ is exceptional.
\end{defn}

\begin{defn}
Let $(X, x)$ be a normal singularity and $\varphi\colon Y\to X$ a
blow-up such that the exceptional locus of $\varphi$ contains only
one irreducible divisor, say $S$, and $x\in \varphi(S)$.
Then $\varphi$ is called a \textsl{plt blow-up} (resp. \textsl{lc
blow-up}) of $(X, x)$, if $(Y,S)$ is plt (resp. lc) and
$-(K_{Y}+S)$ is $\varphi$-ample.
\end{defn}

\begin{prop}[{\cite[Theorem 4.9]{prokhorov}}]
\label{prokhorov}
Let $(X,x)$ be a klt singularity and let $\varphi\colon (Y,S) \to
X$ be a plt blow-up of $(X,x)$.
Then the following are equivalent:
\begin{enumerate}
\item
$(X,x)$ is non-exceptional;
\item
there is a boundary $\Xi \geq \dif_{S}(0)$ such that
$N(K_{S}+\Xi)$ is nef (in particular $\bQ$-Cartier) and $(S, \Xi)$
is not klt;
\item
(in dimension $3$ only) there is a regular (i.e. $1$, $2$, $3$,
$4$ or $6$) complement of $K_{S}+\dif_{S}(0)$ which is not klt.
\end{enumerate}
\end{prop}

About the definition of a complement, the reader is asked to refer
to \ref{ldef}.

\begin{prop}
\label{uniqueness of lc blow-up}
Let $(X,x)$ be a log-canonical singularity. If there are two lc
blow-ups which are not isomorphic over $X$, then $(X,x)$ is not
exceptional.
\end{prop}

\begin{proof}[Proof (cf. {\cite[2.7]{m-p2}})]
Assume that $(X,x)$ is exceptional. Let $\varphi\colon (Y,S)\to X$
be a lc blow-up. Since $-(K_Y+S)$ is $\varphi$-ample, the linear
system $-n(K_Y+S)$ is base point free over $X$ for $n\gg 0$. Let
$H\in |-n(K_Y+S)|$ be a general member and let $B:=\frac1n
\varphi(H)$. By Bertini Theorem, $K_Y+S+\frac1nH$ is lc. Since
$K_Y+S+\frac1nH$ is $\bQ$-linearly trivial, $K_X+B$ is lc and
$a(S^{(i)},B)=-1$ for any component $S^{(i)}$ of $S$. Hence $S$ is
irreducible and $a(S,B)=-1$. If $\varphi'\colon (Y',S')\to X$ is
another lc blow-up, then similarly $S'$ is irreducible and we have
a boundary $B'$ such that $K_X+B'$ is lc and $a(S',B')=-1$. We
claim that $S$ and $S'$ define different discrete valuations of
the function field $K(X)$. Indeed, otherwise $\chi\colon
Y\dasharrow Y'$ is an isomorphism in codimension one and
$\chi_*(K_{Y}+S)=K_{Y'}+S'$. Since both $-(K_{Y}+S)$ and
$-(K_{Y'}+S')$ are ample over $X$, we have $Y\simeq Y'$. We may
assume also  that $B$ and $B'$ have no common components.

For $0\le \alpha\le 1$, define the linear function
$\varsigma(\alpha)$ by
\[
a(S,\alpha B+\varsigma(\alpha) B')=
a(S,0)-\alpha\mult_S(B)-\varsigma(\alpha)\mult_S(B')=-1
\]
and put $B(\alpha):=\alpha B+\varsigma(\alpha) B'$. Clearly,
$\varsigma(1)=0$ and $B(1)=B$ (because $a(S,B)=-1$). 
We claim that $K_X+B(\alpha)$ is lc
for all $0\le \alpha\le 1$. Assume the opposite. Then
\[
0<\alpha_0:=\inf_{\alpha\in[0,1]} \{\alpha \mid
K_X+B(\alpha)\quad\text{is lc}\}.
\]
Fix some log resolution of $(X, B+B')$ factoring through $Y$ and
let $E_1,\dots,E_m$ be the new exceptional divisors. The value
$\alpha_0$ can be computed from a finite number of linear
inequalities $a(E_i, B(\alpha))\ge -1$. Therefore $\alpha_0$ is
rational and $K_X+B(\alpha_0)$ is lc. Since $K_X+B(\alpha_0-\ep)$
is not lc for any $\ep>0$, some inequality $a(E_i,
B(\alpha))\ge -1$ is an equality (see \cite[3.12]{kollar}).
Hence $a(E_i, B(\alpha_0))=a(S, B(\alpha_0))=-1$. This contradicts the
exceptionality of $(X,x)$. Thus $K_X+B(\alpha)$ is lc for all
$0\le \alpha\le 1$. In particular, $K_X+\varsigma(0)B'$ is lc.
Since $(X,x)$ is exceptional, $a(S,B')>-1$ and therefore
$\varsigma(0)>1$. On the other hand, $a(S',B')=-1$ and
$a(S',\varsigma(0)B')<-1$, a contradiction.
\end{proof}

\begin{say}
We make use of toric geometry, and terminologies in \cite{fulton}
are used here. Let $N$ be a free abelian group $\bZ^n$ and $M$ its
dual $Hom_{\bZ}(N, \bZ)$. Denote $N\otimes_{\bZ}\bR$ and
$M\otimes_{\bZ}\bR$ by $N_{\bR}$ and $M_{\bR}$ respectively.
We have a canonical pairing $\langle\,\ \rangle \colon
N_{\bR}\times M_{\bR}\to \bR$. Let $\sigma$ be the positive
quadrant $\bR_{\geq 0}^n$ of $N_{\bR}= \bR^n$ and $\sigma^{\vee}$
its dual. Then $\bCn$ is the toric variety corresponding to the
cone $ \sigma$. For a fan $\D$ in $N$, the corresponding toric
variety is denoted by $T_{N}(\D)$.
For a primitive element $\bp\in N$ of a $1$-dimensional cone
$\tau=\bR_{\geq 0}\bp$ in $\D$, the closure
$\overline{orb(\bR_{\geq 0}\bp)}$ is denoted by $\Dp$, which is a
divisor on $T_{N}(\D)$.
\end{say}

\begin{defn}
A monomial $x_{1}^{m_{1}}\cdots x_{n}^{m_{n}}\in
\bC[[x_{1},x_{2},\ldots, x_{n}]]$ is denoted by $x^m$, where
$m=(m_{1},\cdots,m_{n})\in \bZ^n =M$. For a power series
$f=\sum_{m}a_{m}x^m \in \bC[[x_{1},x_{2},\ldots, x_{n}]]$, we
write $x^m\in f$, if $a_{m}\neq 0$.
For $\bp\in N_{\bR}$ and a power series $f$, we define \[ \bp(f) =
\min_{x^m\in f}\langle\bp,m\rangle. \]
We denote the leading term
$\sum_{\langle\bp,m\rangle=\bp(f)}a_{m}x^{m}$ of $f$ with respect
to $\bp$ by $f_{\bp}$.
\end{defn}

\begin{defn}
For a power series $f=\sum_{m}a_{m}x^m \in
\bC[[x_{1},x_{2},\ldots, x_{n}]]$, define the \textsl{Newton
polyhedron} $\Gamma_{+}(f)$ in $M_{\bR} $ as follows:
\[ \Gamma_{+}(f) =
\text{the convex hull of}\ \bigcup_{x^m\in f}(m+\sigma^{\vee}). \]
The set of the interior points of $\Gamma_{+}(f)$ is denoted by
$\Gamma_{+}(f)^0$.
For each face $\gamma$ of $\Gamma_{+}(f)$, we define the
polynomial $f_{\gamma}$ as follows:
\[ f_{\gamma}=\sum_{m\in\gamma}a_{m}x^m. \]
A power series $f$ is said to be \textsl{non-degenerate}, if for
every face $\gamma$ the equation $f_{\gamma}=0$ defines a
hypersurface smooth in the complement of the hypersurface
$x_{1}\cdots x_{n}=0$.
\end{defn}

\begin{prop}[{\cite[10.3]{varchenko}}]
Let $(X,0)\subset (\bCn, 0)$ be a singularity defined by a
non-degenerate power series $f$. Then there exists a subdivision
$\D$ of $\sigma$ in $N$ such that the toric morphism
$\varphi\colon T_{N}(\D)\to \bCn$ gives an embedded resolution of
$(\bCn, X)$ such that the exceptional set is of pure codimension
one and the union of the proper transform of $X$ and the
exceptional divisor is of normal crossings.

\end{prop}
Here we state a well-known criterion for a non-degenerate
hypersurface singularity to be canonical.

\begin{prop}[see for example, \cite{w-h} or {\cite[Theorem~3]{markushevich}}
or {\cite[Corollary 1.7]{ishii2}}]
\label{markushevich}
Let $(X,0)\subset ( \bCn,0)$ be a normal hypersurface singularity
defined by a non-degenerate power series $f$ and $\Gamma_{+}(f)$
its Newton polyhedron.
Then $(X,0)$ is canonical (resp. log-canonical or lc) if and only
if $\bl=(1, 1, \ldots, 1)\in \Gamma_{+}(f)^0$ (resp. $\bl=(1, 1,
\ldots, 1)\in \Gamma_{+}(f)$), which is equivalent to that
$\bq(f)<\langle \bq, \bl \rangle$ (resp. $\bq(f)\leq \langle \bq,
\bl \rangle$ ) for all $\bq\in \bR^n_{\geq 0}$.

\end{prop}

\section{Plt blow-ups of hypersurface canonical singularities}

\begin{say}
\label{notation}
Let $(X,0) \subset (\bCn, 0)$ be a hypersurface singularity and
$\bp=(p_{1},\ldots p_n)$ a primitive element in $N$ with $p_{i}>0$
for all $i$ (such $\bp$ is called a \textsl{weight}).
Let $\varphi\colon \bCnp\to \bCn$ be the weighted blow-up with a
weight $\bp$. Denote the proper transform of $X$ on $\bCnp$ by
$\Xp$.
The weighted blow-up $\varphi\colon \bCnp\to \bCn$ and its
restriction $\Xp\to X$ are sometimes called the $\bp$-blow-up.
The $\bp$-blow-up $\varphi\colon \bCnp\to \bCn$ is obtained by a
subdivision of the cone $\sigma$. The corresponding fan consists
of the faces of cones $\sigma_{i}$ $(i=1,\ldots,n)$, where
$\sigma_{i}$ is generated by $
\be_{1},\ldots,\be_{i-1},\bp,\be_{i+1},\ldots,\be_{n}$.
Here $\be_{1},\ldots, \be_{n}$ are the unit vectors
$(1,0,\ldots,0),\ldots,(0,\ldots,0,1)$ which generate $\sigma$.
\end{say}

\begin{lem}
\label{lc}
Under the notation of \ref{notation}, let $(X,0)$ be canonical and
defined by a non-degenerate power series $f$.
\begin{enumerate}
\item
If $\bl\in \Gamma_{+}(\fp)$, then $(\bCnp, \Xp+\Dp)$ is lc and
$\Xp$ is normal.
\item
If moreover $\bl\in \Gamma_{+}(\fp)^0$, then $(\bCnp, \Xp+\Dp)$ is
dlt and $\dif_{\Xp}(0)=0$.
\end{enumerate}
\end{lem}

\begin{proof}
Let $\psi\colon \tnd\to \bCnp$ be a toric morphism which is a
log-resolution of $(\bCnp, \Xp+\Dp)$.
Denote
\[ K_{\tnd}+X'+\Dp=\psi^*(K_{\bCnp}+\Xp +\Dp)+\sum_{\bR_{\geq 0}\bq\in
\D(1),\bq\neq \bp}\alpha_{\bq}\Dq,\]
where $X'$ is the proper transform of $X$ on $\tnd$ and by abuse
of notation the divisors corresponding to $\bp$ on $ \bCnp$ and on
$\tnd$ are both denoted by $\Dp$. Here we may assume that $\bq$
are all primitive.
Then if $\bq\in \sigma_{i}$, we obtain that
\[
\alpha_{\bq}=\langle\bq,\bl\rangle-\bq(f)-
\frac{q_{i}}{p_{i}}(\langle\bp,\bl\rangle-\bp(f))-1.
\]
This is because the discrepancy of $K_{\bCnp}+\Xp$ at $\Dq$ is
\[
a( \Dq, \Xp) = \langle\bq,\bl\rangle-\bq(f)-1-
\frac{q_{i}}{p_{i}}(\langle\bp,\bl\rangle-\bp(f)-1)
\]
by \cite[2.6]{ishii2} and the coefficient of $\Dq$ in $\psi^*\Dp$
is $q_{i}/p_{i}$.
Now consider the discrepancy $\alpha_{\bq}$. Since
$\langle\bq,\bl\rangle-\bq(f)>0$ by Proposition
\ref{markushevich}, it follows that $\alpha_{\bq}\geq 0$, if
$q_{i}=0$.
In case $q_{i}>0$, define $\ba:=(a_{1},\ldots,a_{n})$ from the
proportion
$(q_{1}:\ldots:q_{n})=(p_{1}+a_{1}:\ldots:p_{n}+a_{n})$. Then
$a_{i}=0$ and $a_{j}\geq 0$ for $j\neq i$.
Since $\bq(f)\leq \bq(\fp)$, it follows that
\[ \frac{\langle\bq,\bl\rangle-\bq(f)}{q_{i}}\geq
\frac{\langle\bq,\bl\rangle-\bq(\fp)}{q_{i}}= \frac
{\langle\bp+\ba,\bl\rangle-(\bp+\ba)(\fp)}{p_{i}}. \]
Here if one assumes that $\bl\in \Gamma_{+}(\fp)$, the right hand
side of the inequality above is
\[ \frac {\langle\bp,\bl\rangle-\bp(\fp)}{p_{i}}+
\frac {\langle\ba,\bl\rangle-\ba(\fp)}{p_{i}}\geq \frac
{\langle\bp,\bl\rangle-\bp(\fp)}{p_{i}}= \frac
{\langle\bp,\bl\rangle-\bp(f)}{p_{i}},\] which yields that
$\alpha_{\bq}\geq -1$.
If one assumes that $\bl\in \Gamma_{+}(\fp)^0$, then the
inequality above is strict, since
$\langle\ba,\bl\rangle-\ba(\fp)>0$.
So the first assertions of (i) and (ii) are proved.
For the second assertion of (i), consider the discrepancy
$a(\Dq,\Xp)$. This value is $\alpha_{\bq}+q_{i}/p_{i}$ which is
greater than $-1$, if $ q_{i}\neq 0$ and non-negative if $q_{i}=0$
by the argument above.
Hence $(\bCnp, \Xp )$ is plt.
Then, by the inversion of adjunction (\cite[17.6]{uta} or
\cite[2.7]{ishii2}), $\Xp$ is normal and $K_{\Xp}+\dif_{\Xp}(0)$
is klt.
For the second assertion of (ii), note that $(\bCnp, \Xp+\Dp)$ is
dlt and $\bCnp$ is smooth outside of $\Dp$.
Then $\bCnp$ is smooth at a general point of $\Xp\cap\Dp$ by the
classification of $2$-dimensional log-canonical pairs
(\cite[9.6]{kawamata}).
Therefore $\bCnp$ is smooth in codimension two along $\Xp $, which
yields $\dif_{\Xp}(0)=0$.
\end{proof}

\begin{prop}
\label{canonical}
Under the notation of \ref{notation}, let $(X,0)$ be canonical and
defined by a non-degenerate power series $f$.
\begin{enumerate}
\item
If $\bl\in \Gamma_{+}(\fp)$ and $\dif_{\Xp}(0)=0$, then the
singularity $(X_{\bp},0)\subset (\bCn, 0)$ defined by $\fp $ is
log-canonical and $\Xp\to X$ is an lc blow-up.
\item
If $\bl\in \Gamma_{+}(\fp)^0$, then the singularity
$(X_{\bp},0)\subset (\bCn, 0)$ is again canonical and $\Xp \to X$
is a plt blow-up.

\end{enumerate}
\end{prop}

\begin{proof}
For the first assertions of (i) and (ii), by Proposition
\ref{markushevich}, it is sufficient to prove that $(X_{\bp},0)$
is smooth in codimension one under the conditions
$\dif_{\Xp}(0)=0$ and $\bl\in \Gamma_{+}(\fp)$.
If $X_{\bp}$ has a singular locus in codimension one, then
$S=\Dp\cap \Xp$ has a singular locus in codimension one.
Indeed, the restriction $\pi'\colon X_{\bp}\setminus \{0\}\to S$
of the canonical projection $\pi\colon \bCn\setminus \{0\}\to
\bP(p_{1},\ldots,p_{n})=\Dp$ is an open map.
Then we obtain that $\pi'$ is a smooth morphism on the smooth
locus of $S$, as $X_{\bp}$ is a Cohen-Macaulay variety and each
fiber of $\pi'$ is smooth.
Here the singular locus of $S$ is contained in the invariant
divisor $(x_{1}\cdots x_{n}=0)\subset \Dp$, because $S$ is defined
by $\fp$ in the weighted projective space $
\Dp=\bP(p_{1},\ldots,p_{n})$ and $f$ is non-degenerate.
Let $\Sigma$ be a $1$-codimensional component of the singular
locus of $S$ and contained in an invariant divisor $
C=\{x_{i}=0\}$ of $\Dp$.
For simplicity, let $i = 1$. Let $\psi\colon \tnd\to \bCnp$ be a
log-resolution of $(\bCnp, \Xp+\Dp)$ which factors through the
blow-up of $C$ and $X'$ the proper transform of $\Xp$.
Since $\dif_{\Xp}(0)=0$, $(\Xp, S)$ is lc by the previous lemma
and the adjunction. From the classification of $2$-dimensional
log-canonical pairs (\cite[9.6]{kawamata}) there exists an
exceptional divisor $\Dq$ mapped onto $C$ such that the
discrepancy $a(\Dq|_{X'}, S)=-1$.
If $\Dq$ is mapped onto $C=\{x_{1}=0\}$, it follows that
$\bq=\lambda\bp+\mu\be_{1}$, where $\lambda$, $\mu$ are positive
rational numbers.
Since $a(\Dq|_{X'}, S)=\alpha_{\bq}=\langle\bq,\bl\rangle-\bq(f)-
({q_{i}}/{p_{i}})(\langle\bp,\bl\rangle-\bp(f))-1 =-1$ for $i\neq
1$, it follows that $\langle \be_{1},\bl
\rangle-\be_{1}(f_{\bp})=0$, which shows that the order of $\fp$
is $1$, a contradiction.
Now the first assertions of (i) and (ii) are proved.
In particular $\fp$ is irreducible, therefore $S$ is irreducible
and reduced.
Noting that $-(K_{\Xp}+S )=(
\bp(f)-\langle\bp,\bl\rangle)\Dp|_{\Xp}$ is $\varphi$-ample, one
obtains that $\Xp\to X$ is a lc blow-up.
If $\bl\in \Gamma_{+}(\fp)^0$, then, by Lemma ~\ref{canonical},
(ii), $(\Xp, S)$ is dlt, therefore it is plt.

\end{proof}

\begin{cor}
\label{plt blow-up}
Let $(X,0)\subset (\bCn, 0)$ be a canonical singularity defined by
a non-degenerate weighted homogeneous polynomial of weight $\bp$.
Then $\bp$-blow-up $\Xp \to X$ gives a plt blow-up of $ (X,0)$.
\end{cor}

\begin{proof}
In this case, $f=\fp$, therefore the condition of (ii) of
Proposition ~\ref{canonical} holds.
\end{proof}

\begin{thm}
\label{correspondence}
Let $(X,0)\subset (\bCn, 0)$ be an arbitrary hypersurface
canonical singularity defined by a non-degenerate power series
$f$. Then
\begin{enumerate}
\item
there exists a weight $\bp$ such that the weighted blow-up $\Xp\to
X$ gives a plt blow-up of $(X,0)$ and
\item
for a weight $\bp$ obtained in {\rm (i)}, a singularity
$(X_{\bp},0)$ defined by a weighted homogeneous polynomial $\fp$
is again canonical and $(X_{\bp},0)$ is exceptional if and only if
$(X,0)$ is exceptional.
\end{enumerate}
\end{thm}

\begin{proof}
For the statement (i), it is sufficient to show the existence of a
compact face $\gamma$ of $\Gamma_{+}(f)$ such that $\bl\in
\Gamma_{+}(f_{\gamma})^0$, because for every compact face $\gamma$
there exists a weight $ \bp$ such that $f_{\gamma}=f_{\bp}$.
Since \[ \Gamma_{+}(f)=\bigcup_{\gamma:\text{compact face of
}\Gamma_{+}(f)}\Gamma_{+}(f_{\gamma}),\] there exists a compact
face $\gamma$ such that $\bl\in \Gamma_{+}(f_{\gamma})$.
Assume that $\bl$ is a boundary point of $\Gamma_{+}(f_{\gamma})$
for every such $\gamma$ as above. Then $\bl$ belongs to the
non-compact face of $\Gamma_{+}(f_{\gamma})$.
So for any such $\gamma$ there exists $i_{\gamma}$ such that $
H_{i_{\gamma}}:=\{(m_{1},\ldots,m_{n})\in M_{\bR}\mid
m_{i_{\gamma}}\geq 1, m_{j}\geq 0\ (j\neq i_{\gamma})\}$ contains
$\Gamma_{+}(f_{\gamma})$ and $\bl$ is on the boundary of
$H_{i_{\gamma}}$.
Since $\bl$ is on the boundary of $\bigcup_{\bl\in
\Gamma_{+}(f_{\gamma}) }H_{i_{\gamma}}$, it is on the boundary of
$\bigcup_{\bl\in \Gamma_{+}(f_{\gamma})} \Gamma_{+}(f_{\gamma})$,
therefore on the boundary of $\bigcup_{\gamma:\text{compact face
of }\Gamma_{+}(f)}\Gamma_{+}(f_{\gamma})$, a contradiction.
For the first statement of (ii), note that one can take a weight
$\bp$ such that $\bl\in \Gamma_{+}(f_{\bp})^0$ by the argument
above and apply Proposition \ref{canonical}. For such $\bp$,
denote the proper transform of $X_{\bp}$ under the $\bp $-blow-up
$\varphi\colon \bCnp\to \bCn$ by $X_{\bp} (\bp)$ and $X_{\bp}
(\bp) \cap \Dp$ by $S_{\bp}$.
Then in $\Dp$, $S$ and $S_{\bp}$ coincide and $ \dif_{S/ \Xp }(0)$
and $\dif_{S_{\bp}/X_{\bp}(\bp)}(0)$ coincide, because the both
are equal to $\dif_{S/\Dp}(0)+\dif_{\Dp/\bCnp}(0)|_{S}$.
Therefore the conditions for $(X,0)$ and $(X_{\bp},0)$ to be
non-exceptional are the same (Proposition \ref{prokhorov}).
\end{proof}

On determining exceptional hypersurface singularities, now we can
reduce the problem into the weighted homogeneous case.

\section{Weights of weighted homogeneous exceptional singularities}

\begin{lem}
\label{unique}
Let $(X,0)\subset (\bCn, 0)$ be a hypersurface singularity defined
by a power series $f$ and let $\bp$, $\bq$ be two weights such
that neither $\fp$ nor $\fq$ is a power of a single coordinate.
If there is an isomorphism $\Xp\simeq \Xq$ over $X$, then
$\bp=\bq$.
\end{lem}

\begin{proof}
Let $\psi\colon \bCnp(\bq)\to \bCnp$ be the toric morphism
corresponding to the star-shaped decomposition by adding a
$1$-dimensional cone $\bR_{\geq 0}\bq$ and $X'$ be the proper
transform of $\Xp$ by $\psi$. First we show that $\Xp\cap\psi(\Dq)
=\emptyset$.
Assume the contrary.
If $\Dp\cap \Xp\not\subset \psi(\Dq)$, then $\Dp|_{X'}$ and
$\Dq|_{X'}$ are $\bQ$-divisors with different supports, a
contradiction to $\Xp\simeq\Xq$.
If $\Dp\cap \Xp\subset \psi(\Dq)$, then both sides coincide,
because the left hand side is of dimension $n-2$, while the right
hand side is irreducible and of dimension $\leq n- 2$.
Therefore the support of $\Dp\cap \Xp$ is an invariant divisor on
$\Dp$, which yields that $\fp$ is a power of a single coordinate,
a contradiction.
Now we obtain $\Xp\cap\psi(\Dq)=\emptyset$.
It implies that the coefficient of $\Dq$ in $\psi^*(\Xp)$ is 0.
Let $\bp=(p_{1},\ldots,p_{n})$ and $\bq=(q_{1},\ldots,q_{n})$.
If we put $q_{1}/p_{1}:=\min_{i}(q_i/p_{i})$, then $\bq\in
\sigma_{1}$ and the coefficient of $\Dq$ in $ \psi^*(\Xp)$ is
$\bq(f)-\bp(f)(q_{1}/p_{1})$ by \cite[2.6, (2)]{ishii2}.
Hence $\bq(f)-\bp(f)(q_{i}/p_{i})\leq 0$ for all $i$.
Now we obtain that $\bq(f)/q_{i}\leq \bp(f)/p_{i}$ for all $i$.
By making the same procedure with exchanging the role of $\bp$ and
$\bq$, we obtain the opposite inequality $\bq(f)/q_{i}\geq
\bp(f)/p_{i}$ for all $i$, which yields $\bp=\bq$.
\end{proof}

For the assertion of the lemma above, the condition of $\fp$ and
$\fq$ is necessary. In fact we have the following example.

\begin{exmp}[by Tomari]
Let $(X,0)\subset (\bC^3,0)$ be a singularity defined by
$f=x_{1}^k+x_{2}^{k+1}+x_{3}^{k+1}=0$ for $k\geq 3$. Take
$\bp=(1,1,1)$, $\bq=(k+1,k,k)$. Then $\fq=f$ is irreducible and
$\fp=x_{1}^k$ is a power of a single coordinate $x_{1}$. And both
$\Xp$, $\Xq$ are isomorphic to the canonical model over $X$.
\end{exmp}

\begin{prop}
For an exceptional canonical singularity $(X,0)\subset (\bCn, 0) $
defined by a non-degenerate power series $f$, a weight which gives
a plt blow-up is unique. Therefore a weight $\bp$ such that $\bl
\in\Gamma_{+}(\fp)^0$ is unique.
A fortiori, a set $\{m \mid x^m\in \fp\}$ spans a hyperplane in
$M_{\bR}\simeq \bR^n$.
\end{prop}

\begin{proof}
By Proposition \ref{uniqueness of lc blow-up}, plt bow-up for
$(X,0)$ is unique.
If a weight $\bp$ gives a plt blow-up, then by (ii) of Theorem
\ref{correspondence}, $\fp$ is not a power of a single coordinate.
Then apply Lemma \ref{unique}.
\end{proof}

To prove the finiteness of the number of weights for exceptional
singularities defined by weighted homogeneous functions, we need
the following lemma:

\begin{lem}
\label{another weight}
Let $(X,0)\subset (\bCn, 0)$ be a canonical singularity defined by
a non-degenerate power series $f$. Let $g$ be an irreducible
weighted homogeneous polynomial with weight $\bp $. Assume that
$\bl\in \Gamma_{+}(\fp)$ and $x^m\in g$ for every $m$ with $x^m\in
\fp$.
Then $\bp$-blow-up $\Xp\to X$ is a lc blow-up.
\end{lem}

\begin{proof}
By (i) in Proposition~ \ref{canonical}, it is sufficient to prove
that $\dif_{\Xp}(0)=0$.
Assume that $\dif_{\Xp}(0)\neq 0$.
Then a $2$-codimensional irreducible component of the singular
locus of $\bCnp$ is contained in $\Xp$, which implies that
$\{x_{i}=0\}$ in $\Dp$ is contained in $ \{\fp=0\}$ in $\Dp$ for
some $i$.
Therefore $x_{i}\mid \fp$. By $\bl\in \Gamma_{+}(\fp)$, it follows
that $x^2_{i} \nmid \fp$. For the simplicity, let $i=1$.
Here $\bCnp$ is singular along $\{x_{1}=0\}$ in $\Dp$, so the
weight $\bp$ is represented by $(p_{1},ap'_{2},\ldots,ap'_{n})$
for some integer $a\geq 2$ with $a\nmid p_{1}$.
Now, by the assumption on $g$, it follows that $g=x_{1}g'+g''$,
where $x_{1}\nmid g'$ and $x_{1}$ does not appear in $g''$, as $g$
is irreducible.
Then $\bp(g)=\bp(x_{1}g')\equiv p_{1}~\mod a$ and
$\bp(g)=\bp(g'')\equiv 0~\mod a$, which is a contradiction.
\end{proof}

\begin{thm}
\label{finiteness}
For a fixed $n$, the number of weights of non-degenerate weighted
homogeneous polynomials which define exceptional canonical
singularities in $(\bCn, 0)$ is finite.
\end{thm}

\begin{proof}
Let ${\mathcal F}=\{f \mid$ a non-degenerate weighted homogeneous
power series defining an exceptional canonical singularity in $
\bCn \text{ at }0 \}$ and $\bW=\{ \text{ the weight of } f \mid
f\in {\mathcal F}\}.$
Assume that the set $\bW$ is infinite and induce a contradiction.
For $m=(m_{1},\ldots,m_{n})\in \bZ^n_{\geq 0}$, define $
|m|:=\sum_{i=1}^nm_{i}$.

\textbf{Step 1.} For $f\in {\mathcal F}$, define
$\alpha_{f}:=\min\{|m|\mid x^m\in f\}$.
Then $\{\alpha_{f}\}_{f\in {\mathcal F}}$ is bounded by $n$ from
above, since $\bl\in \Gamma_{+}(f)^0$.
Hence there exist a subset ${\mathcal F}_{1}\subset {\mathcal F}$
and $m^{(1)}\in \bZ^n_{\geq 0}$ such that $\bW_{1}:=\{ \text{ the
weight of } f \mid f\in {\mathcal F}_{1}\}$ is a infinite set and
$x^{m^{(1)}}\in f$, $|m^{(1)}|=\alpha_{f}$ for any $f\in {\mathcal
F}_{1}$, because $\{m\in \bZ^n_{\geq 0} \mid |m|\leq n\}$ is
finite and $\bW$ is infinite.
Therefore every $f\in {\mathcal F}_{1}$ is written as $
f=a_{1}x^{{m^{(1)}}}+g$, $(x^{{m^{(1)}}}\notin g,\ a_{1}\in \bC)$.

\textbf{Step2.} For $f=a_{1}x^{{m^{(1)}}}+g\in {\mathcal F}_{1}$,
define $\beta_{f}:=\min\{|m| \mid x^m\in g\}$.
If $\{\beta_{f}\}_{f\in {\mathcal F}_{1}}$ is bounded, then, by
taking infinite sets ${\mathcal F}_{2}$ and $\bW_{2}$ smaller, we
can write $f=a_{1}x^{m^{(1)}}+a_{2}x^{m^{(2)}}+g_{2}$ for every
$f\in {\mathcal F}_{2}$ in the same way as in Step 1.
In particular, if $x^{m^{(1)}}$ is a power of a single coordinate,
then $\{\beta_{f}\}$ is bounded.
Indeed if $\{\beta_{f}\}$ is unbounded, then
$\bigcap_{f\in{\mathcal
F}_{1}}\Gamma_{+}(f)=\Gamma_{+}(x^{m^{(1)}})$, where $\bl$ belongs
to the left hand side and not to the right hand side, a
contradiction.

\textbf{Step 3.} By the successive procedures, one obtains
infinite sets ${\mathcal F}_{r}$, $\bW_{r}$ such that for every $
f\in {\mathcal F}_{r}$, $f=a_{1}x^{m^{(1)}}+\ldots
+a_{r}x^{m^{(r)}} +h$ and $\{\gamma_{f}\}_{f\in{\mathcal F}_{r}}$
is unbounded, where $\gamma_{f}:=\min\{|m|\mid x^m\in h\}$.
Indeed, if these procedures do not terminate, one obtains an
infinite series $\{ m^{(1)}, m^{(2)},\ldots, m^{(i)},\ldots \}$ in
$\bZ^n_{\geq 0}\subset M_{\bR}$.
Let $L_{i}$ be the linear subvariety spanned by $\{ m^{(1)},
m^{(2)},\ldots, m^{(i)}\}$.
Then $L_{i}\neq M_{\bR}$ for all $i$, because $x^{ m^{(1)}},
x^{m^{(2)}},\ldots, x^{m^{(i)}}$ appear in the weighted
homogeneous polynomials in ${\mathcal F}_{i}$.
Since $\ldots\subset L_{i}\subset L_{i+1} \subset \ldots$, there
exists $r$ such that $L_{r}=L_{r+i}$ for all $i>0$.
Let $\bp=(p_{1},\ldots,p_{r})$ be the weight of an element $f\in
{\mathcal F}_{r}$.
Then $L_{r}$ is contained in a hyperplane
$H_{\bp}:=\{(m_{1},\ldots,m_{n})\mid \sum p_{i}m_{i}=p\}$ for some
$p\in \bN$.
Hence the infinite set $\{ m^{(1)}, m^{(2)},\ldots, m^{(i)},\ldots
\}$ is contained in $H_{\bp}\cap \bZ^n_{\geq 0}$ which is a finite
set because of $p_{i}>0$, (i=1,\ldots,n), a contradiction.

\textbf{Step 4.} Now fix an element $f'\in {\mathcal F}_{r}$ and
let $\bp$ be a weight of $f'$.
As $\{\gamma_{f}\}_{f\in{\mathcal F}_{r}}$ is unbounded, one can
take $f=a_{1}x^{m^{(1)}}+\ldots +a_{r}x^{m^{(r)}} +h\in {\mathcal
F}_{r}$ such that $\bp(h)>\bp(f')=\bp(x^{m^{(i)}})$ for
$i=1,\ldots,r$.
Then $\fp=a_{1}x^{m^{(1)}}+\ldots +a_{r}x^{m^{(r)}}$.
So if $x^{m}\in \fp$, then $x^{m}\in f'$.
On the other hand, from the unboundedness of $\gamma_{f}$, it
follows that $\bl \in \bigcap_{f''\in {\mathcal
F}_{r}}\Gamma_{+}(f'') = \Gamma_{+}(\fp)$.
Now by Lemma \ref{another weight}, $\bp$-blow-up of $X=\{f=0\}$ is
an lc blow-up.
On the other hand, for the weight $\bq$ of $f$, $\bq$-blow-up of
$X$ is a plt blow-up by Corollary \ref{plt blow-up}.
Since neither $\fp$ nor $\fq$ is a power of a single coordinate,
by Lemma~\ref{unique} $\Xp\not\simeq \Xq$, which is a
contradiction to that $(X,0)$ is exceptional.
\end{proof}

\begin{rem}
The theorem is the same as the finiteness of Newton polyhedrons of
such singularities.
% If one drops the condition of weighted homogeneous,
% then the number of Newton polyhedrons of such singularities
% is not finite as is seen in
% the following example.
\end{rem}

\begin{cor}
Fix $n\in\bN$. Let $(X,0)\subset(\bCn, 0)$ be a canonical
exceptional singularity defined by a non-degenerate power series
$f$ and let $\varphi\colon (Y,S)\to X$ be a plt blow-up. Then the
pair $(S,\dif_S(0))$ contained in a finite number of algebraic
families.
\end{cor}
Note that this fact is known to be true also for any klt
singularity of dimension $\le 3$ (see \cite[\S 4]{shokurov}).
\begin{proof}
By Theorems~\ref{correspondence} and \ref{finiteness} we may
assume (up to finite numbers of cases) that $\varphi$ is a
weighted blow-up of fixed weight $\bp$. Then the exceptional
divisor $S$ is defined in the weighted projective space $\bP(\bp)$
by $f_{\bp}=0$. Thus we may assume that $S$ is contained in some
algebraic family. Now let $H$ be a very ample divisor on $S$.
Write $\dif_S(0)=\sum (1-1/m_i)\D_i$, where $m_i\in\bN$ and
$\D_i$'s are prime divisors. Since $-(K_S+\dif_S(0))$ is ample,
$\operatorname{Const}\ge H^{n-3}\cdot
(-K_S)>\sum(1-1/m_i)H^{n-3}~\cdot ~\D_i \ge\frac12 \sum
H^{n-3}\cdot \D_i$. Thus the degree of components of $\dif_S(0)$
under the embedding $S\hookrightarrow \bP^N$ given by $|H|$ is
bounded. Then $\Supp(\dif_S(0))$ belongs to a finite number of
families (see, e.g., \cite[Ch. 3 \S 7]{Danilov}) and we may assume
that $(S, \Supp(\dif_S(0)))$ is fixed. Now we need to show only
that $m_i\le \operatorname{Const}$ for all $i$. Indeed, in the
opposite case we can take an infinite sequence $\dif_S^{(1)}(0)<
\dif_S^{(2)}(0)<\cdots$. Let $\dif_S^{(\infty)}(0):=\lim
\dif_S^{(k)}(0)$. Then $-(K_S+\dif_S^{(\infty)}(0))$ is nef and
$\down{\dif_S^{(\infty)}(0)}\ne 0$. This contradicts
Proposition~\ref{prokhorov}.
\end{proof}

\section{Exceptional canonical singularities of Brieskorn type}
The aim of this section is to prove the following:
\begin{thm}
\label{Brieskorn}
Let $X\subset\bC^4$ be a hypersurface canonical singularity given
by the equation
\begin{equation}
\label{eq-first}
x_1^{a_1}+x_2^{a_2}+x_3^{a_3}+x_4^{a_4}=0,\quad a_1\le a_2\le
a_3\le a_4.
\end{equation}
Then $(X,0)$ is exceptional if and only if
$\x{a_1}{a_2}{a_3}{a_4}$ is one of the following:
\par
\begin{tabular}{llp{1cm}ll}
$[3,3,4,d],$&$5\le d\le 11$&&$[3,3,5,d],$&$d=6,7$\\
$[3,4,4,d],$&$d=4,5$&&$[2,3,7,d],$&$8\le d\le 41$\\
$[2,3,8,d],$&$8\le d\le 23$&&$[2,3,9,d],$&$9\le d\le 17$\\
$[2,3,10,d],$&$10\le d\le 14$&&$[2,3,11,d],$&$d=11,12,13$\\
$[2,4,5,d],$&$6\le d\le 19$&&$[2,4,6,d],$&$6\le d\le 11$\\
$[2,4,7,d],$&$d=7,8,9$&&$[2,5,5,d],$&$5\le d\le 9$\\
$[2,5,6,d],$&$d=6,7$.&&&\\
\end{tabular}
\par
Moreover, exceptional divisors $(S,\dif_S(0))$ of corresponding
plt blow-ups are described in Tables~\ref{table-1} and
\ref{table-2}.
\end{thm}

From Proposition~\ref{markushevich} (see also \cite{Reid},
\cite[8.14]{kollar}) we have

\begin{cor}
\label{Bri-can}
Let $X\subset \bC^n$ be a hypersurface singularity
\begin{equation}
\label{eq-Brieskorn}
x_1^{a_1}+x_2^{a_2}+\cdots +x_n^{a_n}=0.
\end{equation}
Then it is canonical (resp. log-canonical) if and only if
\[
1/a_1+1/a_2+\cdots +1/a_n>1,\qquad \text{resp.} \qquad\ge 1.
\]
Moreover, $K_{\bC^n}+X$ is lc if and only if $1/a_1+\cdots
+1/a_n\ge 1$.
\end{cor}

\begin{lem}
\label{MP}
Let $X\subset\bC^n$ be a canonical hypersurface singularity
\eref{eq-Brieskorn}, let $W\subset\bC^n$ be the hyperplane
$\{x_n=0\}$ and let $F:=W\cap X$. If $1/a_1+\cdots +1/a_{n-1}\ge
1$, then $K_X+F$ is lc.
\end{lem}
\begin{proof}
By Corollary~\ref{Bri-can}, $K_W+F$ is lc. Note that $K_{\bC^n}+X$
is plt (by Inversion of Adjunction). Applying \cite[17.7]{uta}
twice we obtain that both $K_{\bC^n}+X+W$ and $K_X+F$ are lc.
\end{proof}

\subsection{Notation}
\label{setup}
From now on we assume that $X\subset \bC^4$ is a hypersurface
canonical singularity given by the equation \eref{eq-first}. The
sequence $\x{a_1}{a_2}{a_3}{a_4}$ is called the \textit{type} of
$X$. Set $w:=\lcm(a_1,\dots,a_4)$ and consider the weighted
blow-up $f\colon Y=X(\bp)\to X$, where $p_i:=w/a_i$. Since $f$ is
weighted homogeneous with respect to $\bp$, $f\colon Y\to X$ is a
plt blow-up (see Corollary~\ref{plt blow-up}). Let $S$ be the
exceptional divisor. Put $\De:=\dif_S(0)$. Note that $S$ is given
in the weighted projective space $\bP(\bp)$ by the equation
\eref{eq-first}. Let $\Ga_i$, $i=1,\dots,4$ be a curve on
$S\subset\bP(\bp)$ which is cut out by $x_i=0$. Set
$\Ga:=\Ga_1+\cdots+\Ga_4$.

\begin{lem}[see e.g. {\cite[8.16]{kollar}}]
\label{bound}
For any $n\in\bN$ there is a constant $\delta(n)$ such that for
any $a_1,\dots, a_n\in\bN$ only one of the following inequalities
holds:
\[
\sum_{i=1}^n\frac1{a_i}\le 1-\delta(n) \qquad\text{or}\qquad
\sum_{i=1}^n\frac1{a_i}\ge 1.
\]
Moreover, $\delta(1)=1/2$, $\delta(2)=1/6$, $\delta(3)=1/42$.
\end{lem}

\subsubsection{}
\label{cases}
If $(X,0)$ is exceptional, then by \cite[Lemma~1.7]{m-p}, $K_X+W$
is not lc for any Cartier divisor $W$ with $W\ni 0$. Thus
\ref{Bri-can} and Lemma~\ref{MP} give us $1/a_1+\cdots+1/a_4>1$
and $1/a_1+1/a_2+1/a_3<1$. Since $(X,0)$ is canonical, $a_1=2$ or
$3$ (see \ref{markushevich}). Then $1/a_2+1/a_3+1/a_4>1/2$ and
$a_2\le 5$. Further, $a_2\ge 3$ and $1/a_3+1/a_4>1/6$. Hence,
$a_3\le 11$. Finally, by Lemma~\ref{bound} we have
$1/a_1+1/a_2+1/a_3\le 41/42$ and $a_4\le 41$. This yields for
$\x{a_1}{a_2}{a_3}{a_4}$ cases as in Theorem~\ref{Brieskorn} and
additionally the following cases: $\x2377$, $\x3344$, $\x2455$. In
these cases the singularity is not exceptional. This will be
proved in \ref{non-exc-2-3-7-7}, \ref{non-exc-3-3-4-4} and
\ref{non-except-2-4-5-5}, respectively.

\begin{remark*}
As above, Lemma~\ref{bound} gives an effective bound of canonical
exceptional singularities of the Brieskorn type in any dimension.
\end{remark*}

\begin{lem}
\label{properties-Gamma}
Let $S\subset\bP=\bP(p_1,\dots,p_n)$ be a hypersurface
\begin{equation*}
%\label{eq-first-1}
x_1^{a_1}+\cdots+x_n^{a_n}=0,\quad\text{where}\quad
a_1p_1=\cdots=a_np_n=w.
\end{equation*}
Let $\Ga_i:=S\cap\{x_i=0\}$, $i=1,\dots,n$ and $\Ga=\sum\Ga_i$.
Then
\begin{enumerate}
\item
$\dif_{S/\bP}(0)=0$;
\item
$K_S+\Ga_i$ is plt for all $i=1,\dots,n$;
\item
$K_S+\Ga$ is lc.
\end{enumerate}
\end{lem}
\begin{proof}
(i) is obvious. Indeed, codimension two singularities of $\bP$ are
contained in $\cup_{i\ne j}\{x_i=x_j=0\}$ and $S$ does not contain
its components. To prove (ii) and (iii) we consider the finite map
$\bP^{n-1}\to \bP$ given by
\[
(x_1,\dots,x_n)\to (x_1^{p_1},\dots,x_n^{p_n}).
\]
The ramification divisor is $\sum(p_i-1)H_i$, where $H_i$ is
$i$-th coordinate hyperplane on $\bP^{n-1}$. Let
$S'=\{x_1^w+\cdots+x_n^w=0\} \subset \bP^{n-1}$ be the preimage of
$S$. The restriction $\var\colon S'\to S$ is also a finite
morphism of degree $w$. Put $L_i:=H_i\cap S'$. By the ramification
formula we have
\[
\var^*\left(K_S+\sum\Ga_i\right)=K_{S'}+\sum L_i.
\]
Since $S'+\sum L_i$ is a normal crossing divisor, $K_{S'}+\sum
L_i$ is lc. Then $K_S+\sum\Ga_i$ is lc by \cite[\S 2]{shokurov2}
or \cite[20.3]{uta}. (ii) can be proved in a similar way.
\end{proof}

\begin{corollary1}
\label{support}
Notation as in \ref{setup}. Let $\Theta$ be a boundary on $S$.
Assume that $\Supp(\Theta)\subset \Ga$ and $\down{\Theta}=0$. Then
$K_S+\Theta$ is klt.
\end{corollary1}

\begin{lem}
\label{different}
Notation as in \ref{setup}. Then $\De=\sum (1-1/m_i)\Ga_i$, where
\[
m_i=\gcd(p_1,\dots,\hat{p_i},\dots,p_4).
\]
\end{lem}
\begin{proof}
Restricting $K_{\bC^4(\bp)}$ to $S$ we obtain
\[
\dif_{S/Y}(0)+\dif_{Y/\bC^4(\bp)}(0)|_S=
\dif_{S/\bP}(0)+\dif_{\bP/\bC^4(\bp)}(0)|_S.
\]
By Lemma~\ref{lc} and Lemma~\ref{properties-Gamma},
$\dif_{Y/\bC^4(\bp)}(0)=0$ and $\dif_{S/\bP}(0)=0$. Hence,
$\dif_{S/Y}(0)=\dif_{\bP/\bC^4(\bp)}(0)|_S$. Taking into account
that $\dif_{\bP/\bC^4(\bp)}(0)=\sum (1-1/m_i)\{x_i=0\}$, we get
the assertion.
\end{proof}

We say that expression $\bp=(p_1,\dots,p_n)$ is
\textit{normalized} if for each $i$ we have
$\gcd(p_1,\dots,\hat{p_i},\dots,p_n)=1$.

\begin{lem}[see, e.g., {\cite[1.3.1]{Dol}}]
\label{normalized}
\begin{enumerate}
\item
If $q=\gcd(p_1,\dots,p_n)$, then
$\bP(p_1,\dots,p_n)\simeq\bP(p_1/q,\dots,p_n/q)$.
\item
If $\gcd(p_1,\dots,p_n)=1$ and $q=\gcd(p_2,\dots,p_n)$, then the
map
\[
\bC^{n+1}\longleftarrow \bC^{n+1},\qquad
(x_1^q,x_2,\dots,x_n)\longleftarrow (\bar x_1,\bar x_2,\dots,\bar
x_n).
\]
induces the isomorphism
\begin{equation*}
\bP(p_1,p_2,\dots,p_n)\simeq\bP(p_1,p_2/q,\dots,p_n/q).
\end{equation*}
\end{enumerate}
\end{lem}

\begin{lem}
\label{canonical-divisor}
Let $\bP=\bP(\bp)$ be a weighted projective space, where
$\bp=(p_1,\dots,p_n)$. Assume that $\bp$ is normalized. Then
\begin{enumerate}
\item
$K_{\bP}\sim \OOO_{\bP}\left(-\sum p_i\right)$;
\item
$\OOO_{\bP}(1)^{n-1}=1/p_1\cdots p_n$.
\end{enumerate}
\end{lem}
\begin{proof}
(i) follows, for example, from the discussion 2.1 and 2.2 of
\cite{Dol}. (ii) follows easily from the fact that $\bP$ is a
finite abelian quotient of the projective space.
\end{proof}

\subsection{}
\label{algorithm}
By Lemma~\ref{normalized} we may assume that the exceptional
divisor $S$ is given by the equation $y_1^{\bar
a_1}+\cdots+y_4^{\bar a_4}=0$ in $\bP(\bar{\bp})$, where
$\bar{\bp}=(\bar p_1,\dots,\bar p_4)$ is normalized. The algorithm
of computation of $\bp$ and $(\bar a_1,\dots,\bar a_4)$ is as
follows.
\par
Starting with $\x{a_1}{a_2}{a_3}{a_4}$ we find
$w=\lcm(a_1,\dots,a_4)$, $p_i=w/a_i$. Then
$\gcd(p_1,\dots,p_4)=1$. For convenience we put
$a_1,\dots,a_4,p_1,\dots,p_4$ into a $(2\times 4)$-matrix and
perform the following transformations: $$
\begin{pmatrix}
a_1&a_2&a_3&a_4\\ p_1'&p_2'&p_3'&p_4'\\
\end{pmatrix}
\longrightarrow
\begin{pmatrix}
a_1/d'&a_2&a_3&a_4\\ p_1'&p_2'/d'&p_3'/d'&p_4'/d'\\
\end{pmatrix},
\leqno 1) $$ where $d':=\gcd(p_2',p_3',p_4')$,\quad 2)\dots,
4)\dots In four steps we get the matrix
\[
\begin{pmatrix}
\bar a_1&\bar a_2&\bar a_3&\bar a_4\\ \bar p_1&\bar p_2&\bar
p_3&\bar p_4\\
\end{pmatrix}
\]
with the normalized second row. Then
\[
S=\{y_1^{\bar a_1}+\cdots+y_4^{\bar a_4}=0\}\subset
\bP(\bar{\bp}).
\]
By Lemma~\ref{different}, $\De=\sum (1-1/m_i)\Ga_i$, where
$m_i=\gcd(p_1,\dots,\hat{p_i},\dots,p_4)$.
\par
Note that if $\bar a_k=1$ for some $k$, then the projection
\begin{equation}
\label{eq-projection}
S\to \bP(\bar p_1,\dots,\hat{\bar p_k},\dots,\bar p_4)
\end{equation}
is an isomorphism. It is easy to see that this holds if and only
if
 $$
\exists k\qquad (a_k,a_i)=1\qquad\text{for all}\qquad i\ne k.
\leqno(\star)
 $$
The following lemma can be easily proved by direct local
computations.
\begin{lem}[cf. {\cite{Di}}]
\label{sing}
Let $S\subset\bP(\bp)$ be a hypersurface \eref{eq-first}, where
$\bp=(p_1, p_2, p_3, p_4)$, $p_1a_1=\cdots=p_4a_4=w$ and $\bp$ is
a weight. Then $\Sing(S)=S\cap \Sing(\bP)\subset \cup_{i\ne
j}\Ga_i\cap\Ga_j$. If $\bp$ is normalized, then at points
$\Ga_i\cap\Ga_j$, $i\ne j$ the surface $S$ has a singularity of
type $\frac1d(p_i,p_j)$, where $d=\gcd(p_k,p_l)$,
$\{k,l\}\cap\{i,j\}=\emptyset$. Moreover,
$\#\Ga_i\cap\Ga_j=w/\lcm(p_k,p_l)$.
\end{lem}

\subsection{Singularities which satisfy $(\star)$}
\label{setup-1}
First consider singularities which satisfy the condition
$(\star)$. Then $S\simeq\bP(\bq)$, where
$\bq=(q_1,q_2,q_3)=(p_1,\dots,\hat{p_k},\dots,p_4)$. Let $L_1$,
$L_2$, $L_3$ be the coordinate lines in $\bP(\bq)$. For $m\in\bN$
let $C_m$ be a curve in $\bP$ given by the equation
$x^{m/q_1}+y^{m/q_2}+z^{m/q_3}=0$ (we assume that $q_i\mid m$ for
$i=1,2,3$). Note that $C_m$ is a smooth curve contained in the
smooth locus of $\bP(\bq)$. The projection \eref{eq-projection}
identifies $\Ga_1,\dots,\hat{\Ga_k},\dots,\Ga_4$ with
$L_1,L_2,L_3$ and $\Ga_k$ with $C_{\bar w}$, where $\bar w=:\bar
a_1\bar p_1=\cdots =\bar a_4\bar p_4$. By $H$ denote the positive
generator of the Weil divisor class group of $S$. If $\bq$ is
normalized, then $\OOO(H)=\OOO(1)$ (see \cite{Dol}). Recall also
that $H^2=1/q_1q_2q_3$. Taking into account \ref{algorithm} and
\ref{cases} we obtain Table~\ref{table-1} and additionally case
$\x2377$ below.

\subsubsection{Case $\x2377$}
\label{non-exc-2-3-7-7}
Then $S=\bP(7,1,1)$ and $\De= \frac12C_7+ \frac23L_1$. Take
$\De^+= \frac12C_7+ \frac23L_1+ \frac56M$, where $M:=\{y+z=0\}$.
Then $K_S+\De^+\qq 0$. It is easy to see that $K_S+\De^+$ is not
klt at $C_7\cap L_1\cap M$. Here the singularity is
non-exceptional.

\subsubsection*{}
Now we prove that all singularities in Table~\ref{table-1} are
exceptional. We consider them case by case according to the type
of the surface $S$. We will assume that there exists a regular
$n$-complement $K_S+\De^+$ and derive a contradiction or prove
that $K_S+\De^+$ is klt (see \ref{prokhorov}). Set
$\De':=\De^+-\De$. We need the definition and a few properties of
complements.

\begin{defn}[\cite{shokurov2}]
\label{ldef}
Let $S$ be a normal variety and let $D=C+B$ be a boundary on $S$,
where $C:=\down{D}$ and $B:=\fr{D}$. Then we say that $K_S+D$ is
\textit{$n$-complementary}, if there is a $\bQ$-divisor $D^+$ such
that
\begin{enumerate}
\item
$n(K_S+D^+)\sim 0$ (in particular, $nD^+$ is an integral divisor);
\item
$K_S+D^+$ is lc;
\item
$nD^+\ge nC+\down{(n+1)B}$.
\end{enumerate}
In this situation an \textit{$n$-complement} of $K_S+D$ is
$K_S+D^+$. We say that an $n$-complement is \textit{regular} if
$n\in\{1,2,3,4,6\}$.
\end{defn}

\begin{lem}[{\cite{shokurov}}, {\cite{Pr_lect}}]
\label{properties-complemnets}
Let $K+\sum\delta_i^+\De_i$ be an $n$-complement of
$K+\sum\delta_i\De_i$, where $\De_i$'s are irreducible components.
\begin{enumerate}
\item
If $\delta_i=1-1/m_i$, $m_i\in\bN$, then $\delta_i^+\ge \delta_i$.
\item
If $\delta_i\ge 6/7$ and $n\in\{1,2,3,4,6\}$, then $\delta_i^+=1$.
\item
If $\delta_i=4/5$ and $n\in\{1,2,3,4,6\}$, then $\delta_i^+=5/6$
or $1$.
\end{enumerate}
\end{lem}

\subsubsection*{Cases when $S\simeq\bP^2$}
For example, assume that $a_2=a_3=a_4$. Here $\gcd(a_1,a_2)=1$.
Algorithm~\ref{algorithm} is as follows:
\[
\begin{pmatrix}
a_1&a_2&a_2&a_2\\ a_2&a_1&a_1&a_1\\
\end{pmatrix}
\longrightarrow
\begin{pmatrix}
1&a_2&a_2&a_2\\ a_2&1&1&1\\
\end{pmatrix}.
\]
Thus $S\simeq\bP^2$ and $\De=\frac{a_1-1}{a_1}C_{a_2}$. There are
two possibilities: $\x3444$ and $\x2555$. In case $\x3444$, among
regular complements of $K_S+\De$ there are only $3$-complement
$K_{\bP^2}+\frac23C_4+\frac13M_1$, where $M_1$ is a line,
$6$-complement $K_{\bP^2}+\frac23C_4+\frac16M_2$, where $M_2$ is a
conic and $6$-complement
$K_{\bP^2}+\frac23C_4+\frac16M_1'+\frac16M_1''$, where $M_1'$,
$M_1''$ are lines. All these complements are klt by
Lemma~\ref{can-term} below. Case $\x2555$ is similar.

\begin{lem}[cf. {\cite[Lemma~3]{prokhorov}}]
\label{can-term}
Let $(S,o)$ be a normal analytic surface germ and let $M=\sum
d_iM_i$ be a $\bQ$-divisor on $S$. Assume that $K_S+M_i$ is plt at
$o$ for all $i$ (for instance, this holds if $(S,o)$ is smooth and
all $M_i$'s also are smooth at $o$). If $\sum d_i\le 1$ and
$\down{M}\le 0$, then $K_S+M$ is klt at $o$.
\end{lem}
\begin{proof}
Let $\var\colon (S',o')\to (S,o)$ be a finite \'etale in
codimension one cover such that $S'$ is smooth. Set $M':=\var^*M$
and $M_i':=\var^*M_i$. By \cite[20.4]{uta}, $K_{S'}+M_i'$ is plt
for all $i$. Hence, all $M_i'$'s are smooth irreducible curves.
Again, by \cite[20.4]{uta} it is sufficient to show that
$K_{S'}+M'$ is klt. By our assumption, $M'=\sum d'_iM_i'$, where
$\sum d_i'\le 1$ and $\down{M'}\le 0$. Let $\sigma\colon\bar{S}\to
S'$ be the blow-up of $o'$ and let $\bar{M}$ be the crepant
pull-back of $M'$ (i.e.
$\sigma^*(K_{S'}+M')=K_{\bar{S}}+\bar{M}$). By
\cite[3.10]{kollar}, $K_{S'}+M'$ is klt if and only if so is
$K_{\bar{S}}+\bar{M}$.  Clearly, all the irreducible components
$\bar{M}_i$ of $\bar{M}$ are smooth. Write $\bar{M}=\sum
\bar{d}_i\bar{M}_i$ so that $\bar{M}_0$ is the exceptional divisor
and $\bar{M}_i$'s are proper transforms of $M'_i$'s for $i\ne 0$.
It is easy to see that $\bar{d}_i=d'_i$ for $i\ne 0$ and
$\bar{d}_0=\sum d'_i-1\le 0$. So we again have $\sum \bar{d}_i\le
1$. Thus, it is sufficient to prove our statement on $\bar{S}$. We
replace $S'$ with $\bar{S}$ and continue the process. At the end,
we get the situation when $\Supp(M')$ is a normal crossing
divisor. In this situation, the inequalities $\sum d_i'\le 1$ and
$d_i'<1$ gives us that $K_{S'}+M'$ is klt (and even canonical).
\end{proof}

If $(a_i,a_j)=1$ for all pairs $(i,j)$, $i\ne j$, then $S\simeq
\bP^2$ and $\De= \sum_{i=1}^3(1-1/a_i)L_i+(1-1/a_4)C_1$. For
$\x{a_1}{a_2}{a_3}{a_4}$ there are the following possibilities:
\[
\x23{11}{13}\quad\text{and}\quad \x237r, r\in\{11, 13, 17, 19, 23,
25, 29, 31, 37, 41\}.
\]
All these singularities are exceptional. Indeed, if $K_S+\De^+$ is
a regular complement of $K_S+\De$, then $\De^+\ge \De$ (see
\ref{properties-complemnets}). Therefore $\De^+\ge
\frac12L_1+\frac23L_3+L_3+C_1$. But $\deg\De^+=3$, a
contradiction. Similarly, one can treat the other cases with
$S\simeq\bP^2$ (see Table~\ref{table-1}).
\begin{remark*}
For any canonical singularity $(X,0)$, define
\[
\compl(X):=\min\{m\mid\ \text{there is a non-klt $m$-complement
of}\ K_X\}.
\]
Let $f\colon (Y,S)\to X$ be a plt blow-up. Then
\[
\compl(X)\le \min\{m\mid K_S+\dif_S(0)\ \text{is
$m$-complementary}\}
\]
(see \cite[Corollary~1]{prokhorov}). Moreover, if $(X,0)$ is
exceptional, then equality holds. It follows from \cite{shokurov}
that $\compl(X)$ is bounded in the three-dimensional case. If $X$
is of the Brieskorn type $\x23{11}{13}$, then $\compl(X)=66$ (see
\cite{Pr_lect}). This is the maximal known value of $\compl(X)$
for three-dimensional canonical singularities. Note that in the
two-dimensional case $\compl(X)\le 6$ and the equality achieves
for singularities of type $E_8$ (= Brieskorn type $[2,3,5]$). By
\ref{prokhorov}, $\compl(X)\in\{1,2,3,4,6\}$ for any
three-dimensional non-exceptional singularity.
\end{remark*}

\begin{conjecture*}
Let $(X,o)$ be a canonical singularity. Then $\compl(X)\le 66$.
\end{conjecture*}
It is known also that the inequality $\compl(X)\le 66$ holds for
any isolated log-canonical three-dimensional singularity
\cite{Is-1}.

\subsubsection*{}
Now we consider cases when $S\ne \bP^2$. In many cases, $S$ is a
cone over a rational normal curve.

\subsubsection*{Cases when $S\simeq\bP(1,1,2)$ (quadratic cone)}
In cases: $\x237{2r}$, $r\in \{4,5,8,10,11,13,16,17,19,20\}$,
$\x238r$, $r\in \{11,13,17,19,23\}$, $\x245r$, $r\in
\{7,9,11,13,17,19\}$, $\x2479$, $\x23{10}r$, $r\in \{11,13\}$,
$\x2567$ and $\x2478$, $\De$ has a component $\De_1\equiv 2H$ with
the coefficient $\ge 6/7$. By Lemma~\ref{properties-complemnets},
for any regular complement $\De^+$ we have $\De^+\ge \De_1$. This
yields a contradiction with $\De^+\equiv-K_S\equiv 2H$.

Similarly, in cases $\x245{4r}$, $r=2,3,4$ and $\x238{20}$, $\De$
has a component $\De_1$ with the coefficient $4/5$. Again
$\De^+\ge\frac56\De_1$. In case $\x245{16}$ this gives us a
contradiction. In cases $\x245{12}$ and $\x238{20}$, we obtain
$\Supp(\De^+)=\Supp(\De)$. By Corollary~\ref{support}, $K_S+\De^+$
is klt. Finally, in case $\x245{8}$ we have only one possibility
with $\Supp(\De^+)\ne\Supp(\De)$:
$\De^+=\frac56C_4+\frac12L_3+\frac16M$, where $M\equiv H$, $M\ne
L_3$. Since $C_4\cap L_3\cap M=\emptyset$, we may apply
Lemma~\ref{can-term}. Here $K_S+M$ is plt by Lemma~\ref{plt-0}
below.

\begin{lem}
\label{plt-0}
Let $(S,o)$ be a germ of a surface singularity of type
$\frac1n(1,1)$ and let $M$ be a germ of a smooth curve passing
through $o$. Then $K_S+M$ is plt at $o$.
\end{lem}
\begin{proof}
Let $\sigma\colon\tilde{S}\to S$ be the minimal resolution, let
$E$ be the (irreducible) exceptional divisor and let $\tilde{M}$
be the proper transform of $M$. Write
$\sigma^*(K_S+M)=K_{\tilde{S}}+\tilde{M}+\alpha E$. Then $\sigma$
is a log-resolution and it is sufficient to show that $\alpha<1$.
By Adjunction, $0=(K_{\tilde{S}}+E)\cdot E+\tilde{M}\cdot
E+(\alpha-1) E^2= -2+1-(\alpha-1)n$. Thus $\alpha=1-1/n$.
\end{proof}

\subsubsection*{Cases when $S\simeq\bP(1,1,3)$ (rational cubic cone)}
Almost all cases can be treated as above because $\De$ has a
component with the coefficient $4/5$ or $\ge 6/7$. The only
non-trivial case is $\x2399$. Then $\De'\equiv\frac12H$, where $H$
is as in \ref{setup-1}. Therefore $C_9$ is not a component of
$\De'$. If $\Supp(\De')\not\subset\sum L_i$, then we have only one
possibility $\De'=\frac16M$, where $M\equiv 3H$, or
$\De'=\sum\alpha_iM_i$, where $M_i$'s are generators of the cone
and $\sum\alpha_i=1/2$. In both cases $K_S+\De^+$ is klt by
Lemma~\ref{can-term}.

\subsubsection*{Case $\x238{8r}$, $r=1,2$}
From $\De_{r=1}\le \De_{r=2}$ we may assume that $r=1$, i.e.,
$\De=\frac23C_8$. Then $\De'\equiv \frac23H$. If the coefficient
of $C_8$ in $\De^+$ is bigger than that in $\De$, then there is
only one possibility $\De^+=\frac34C_8$. Clearly, this complement
is klt. Thus we may assume that $\De^+=\frac23C_8+\De'$, where
$C_8$ is not a component of $\De'$. Note that $S$ is isomorphic to
a projective cone in $\bP^5$ over a rational normal curve of
degree $4$. If all components of $\De'$ are generators, then we
can write $\De'=\sum \alpha_iM_i$, where $\sum \alpha_i=2/3$. By
Lemma~\ref{can-term} and Lemma~\ref{plt-0}, $K_S+\frac23C_8+\De'$
is klt at the vertex in this case. Assume that
$K_S+\frac23C_8+\De'$ is not klt at some point $P$ (outside of the
vertex). Then $P\in C_8$ and there is exactly one component, say
$M_1$ of $\De'$ passing through $P$. By assumption, $C_8$ and
$M_1$ at $P$ cannot intersect each other transversally. Since
$C_8\cdot M_1=8H^2=2$, $C_8$ and $M_1$ have simple tangency at
$P$. Using $\alpha_1\le 2/3$, we can easily check that
$K_S+\frac23C_8+\alpha_1M_1$ is klt at $P$. Finally, assume that
$\De'=\sum \alpha_i M_i$ and $M_1$ is not a generator of the cone.
Then $\alpha_1\ge 1/6$ and $M_1\sim kH$, where $k\ge 4$. This
yields $k=4$ and $\De'=\frac16M_1$. Since $M_1$ is irreducible, it
does not contain the vertex of the cone. By Lemma~\ref{can-term},
$K_S+\frac23C_8+\frac16M_1$ is klt.

\subsubsection*{Cases $\x237{35}$ and $\x3445$}
The only regular complements are $\De^+= \frac12C_7+ \frac23L_1+
\frac56L_3$ and $\De^+=\frac23L_1+ \frac56C_4$, respectively. They
are klt (see Corollary~\ref{support}).

\subsubsection*{Case $\x23{10}{10}$}
Then $\De'\equiv\frac13H$. If the coefficient $\alpha$ of $C_{10}$
in $\De^+$ is bigger than that in $\De$, then $7/10\ge\alpha>2/3$
and $12\alpha\in\bZ$. This is impossible. As in the case $\x2388$
there is only one possibility: $\De^+=\frac23C_{10}+\sum\alpha_i
M_i$, where $\sum \alpha_i=1/3$. By Lemma~\ref{can-term} this
complement is klt.

\subsubsection*{Case $\x23{11}{11}$}
Then $\De'\equiv \frac16H$. The only possibility for regular
complement with $\Supp(\De^+)\ne \Supp(\De)$ is
$\De^+=\frac12C_{11}+\frac23L_1+\frac16M$, where $M:=\{z=\const
x\}$. By Lemma~\ref{plt-0}, $K_S+M$ is plt. We have to check only
that $K_S+\De^+$ is klt near $M$. This follows by
Lemma~\ref{can-term} because $M\cap C_{11}\cap L_1=\emptyset$.

\subsubsection*{Cases when $S\simeq\bP(3,2,1)$ (Gorenstein del Pezzo surface
of degree $6$)} All cases are exceptional because we have a
component $C_6$ of $\De$ with the coefficient $\ge 6/7$.

\subsubsection*{Case $\x237{21}$}
Then $\De'\equiv \frac12H$. Since $K_S+\frac23C_{21}$ is ample,
the coefficient of $C_{21}$ in $\De^+$ is $1/2$. If
$\Supp(\De^+)\not\subset\Ga$, then $\De^+=\frac12C_{21}+\frac16M$,
where $M:=\{x_2=\const x_3^3\}$, $\const\ne 0$. The curve $M$
contains the point $(1,0,0)$ of type $\frac17(1,2)$. Taking into
account that $(K_S+M)\cdot M=-8/7$ and by Lemma~\ref{plt-1} below,
$K_S+M$ is plt. Finally, by Lemma~\ref{can-term}, $K_S+\De^+$ is
klt.

\begin{lem}
\label{plt-1}
Let $S$ be a projective surface with only log terminal
singularities and let $M$ be an irreducible curve on $S$. Assume
that $M$ contains singular points $P_1,\dots,P_r$ of $S$ of types
$\frac1{m_1}(1,q_1),\dots,\frac1{m_r}(1,q_r)$, where
$\gcd(m_i,q_i)=1$, $\forall i$. If
\begin{equation}
\label{ineq}
(K_S+M)\cdot M\le -2+\sum (1-1/m_i),
\end{equation}
then $K_S+M$ is plt near $M$.
\end{lem}
\begin{proof}
Let $\Phi\subset M$ be the set of points where $K_S+M$ is not plt.
Let $\nu\colon\widehat{S}\to S$ be a birational morphism which is
a log resolution over $\Phi$ and an isomorphism outside of $\Phi$.
Write
\[
\nu^*(K_S+M)=K_{\widehat{S}}+\widehat M+\sum a_iE_i,
\]
where $\widehat M$ is the proper transform of $M$ and $\sum
a_iE_i$ is the exceptional divisor. Then $\widehat M$ is smooth
and by Adjunction
\[
\left.\left(K_{\widehat{S}}+\widehat M+\sum
a_iE_i\right)\right|_{\widehat M}= K_{\widehat M}+\dif_{\widehat
M}\left(\sum a_iE_i\right).
\]
We can write
\[
\dif_{\widehat M}\left(\sum a_iE_i\right)=
\sum_{P'\in\nu^{-1}(\Phi)\cap \widehat M} a'P'+
\sum_{P'\notin\nu^{-1}(\Phi)\cap \widehat M} b'P',
\]
If $P\notin\Phi$, then by construction, $\nu$ is an isomorphism
over $P$ and by \cite[16.6]{uta}, $(S,M)\simeq_{an}
(\bC^2,\{x=0\})/\bZ_m(1,q)$, $\gcd(m,q)=1$ and $b'=1-1/m$.

If $P\in\Phi$ then by Connectedness Lemma \cite[5.7]{shokurov2},
$\widehat{M}+\sum_{a_i\ge 1} E_i$ is connected near $\nu^{-1}(P)$.
Thus the coefficient $a'$ of $\dif_{\widehat M}\left(\sum
a_iE_i\right)$ at $P'\in \nu^{-1}(P)$  is $\ge 1$. Now we have
$\dif_{\widehat M}\left(\sum a_iE_i\right)\ge\sum_{i=1}^r
(1-1/m_i)P_i'$, where $\nu(P_i')=P_i$. Combining this with
\eref{ineq} and $\deg K_{\widehat M}\ge -2$, we obtain
$\nu^{-1}(P_i)=\{P_i'\}$ and $\dif_{\widehat M}\left(\sum
a_iE_i\right)=\sum_{i=1}^r (1-1/m_i)P_i'$. In particular,
\[
\down{\dif_{\widehat M}\left(\sum a_iE_i\right)}=0
\]
and $\Phi=\emptyset$.
\end{proof}

\subsubsection*{Cases $\x237{14r}$, $r=1,2$}
Since $\De_{r=1}\le \De_{r=2}$, it is sufficient to consider only
case $r=1$. We have $\De'\equiv \frac23H$. If $C_{14}$ and $L_i$'s
are the only components of $\De^+$, then $K_S+\De^+$ is klt by
\ref{support}. Assume that there is a component $M\ne L_i,
C_{14}$. Then $M\in |rH|$, for some $r$. It is clear that
$-(K_S+\frac12L_3 +\frac23C_{14}+\frac16M)$ is nef. This gives us
$10-\frac23\cdot14-\frac16r\ge 0$ and $r\le 4$. Since $M\ne L_3$,
$r>1$. Further, $M$ contains the point $\{y=z=0\}$ of type
$\frac17(2,1)$. By Lemma~\ref{plt-1}, $K_S+M$ is plt. Then
$K_S+\De^+$ is klt by Lemma~\ref{can-term}.

\subsection{Singularities which do not satisfy $(\star)$}
\label{non-Q-fac}
Now we consider singularities which do not satisfy the condition
of $(\star)$. We will see that all of them are not analytically
$\bQ$-factorial (see \ref{Q-factorial}). Recall that we assumed
that $S\subset\bP(\bar{\bp})$ with normalized $\bar{\bp}$. By $H$
denote the class of Weil divisors such that $\OOO_S(H)=\OOO_S(1)$.
By Lemma~\ref{canonical-divisor},
\[
H^2=\bar w/\bar p_1\cdots\bar p_4, \qquad K_S\sim\left(\bar
w-\sum\bar p_i\right)H,
\]
where $\bar w=\bar a_1\bar p_1=\cdots=\bar a_4\bar p_4$. Taking
into account \ref{algorithm} and \ref{cases} we obtain
Table~\ref{table-2} and additionally cases $\x3344$ and $\x2455$
below.

\subsubsection{Case $\x3344$}
\label{non-exc-3-3-4-4}
Then $S=\{y_1^3+y_2^3+y_3^4+y_4^4=0\}\subset \bP(4,4,3,3)$ and
$\De=0$. We claim that the singularity is not exceptional. Indeed,
consider the curve $M:=S\cap \{y_3=\omega y_4\}$, where
$\omega^4=-1$. We have $K_S+\frac23M\qq 0$. It is sufficient to
show that $K_S+\frac23M$ is not klt (see
Proposition~\ref{prokhorov}). Indeed, $M=\{y_3-\omega
y_4=y_1^3+y_2^3=0\}$ has three components passing through one
point $(0,0,\omega,1)$. This point is singular of type
$\frac13(1,1)$. By blowing-up it we obtain an exceptional divisor
$E$ with discrepancy
\[
a\left(E,\frac23M\right)=a(E,0)-\frac23\Bigl(m_1+m_2+m_3\Bigr)=
-\frac13-\frac23\Bigl(m_1+m_2+m_3\Bigr),
\]
where $m_1$, $m_2$, $m_3$ are multiplicities of components of $M$.
Clearly, $m_1, m_2,m_3\ge 1/3$. Therefore $a(E,M)\le-1$ and
$K_S+\dif_S(\frac23M)$ is not klt.

\subsubsection{Case $\x2455$}
\label{non-except-2-4-5-5}
Then $S=\{\y2255=0\}\subset\bP(5,5,2,2)$ and $\De=\frac12\Ga_2$.
As in case $\x3344$ we can take $\De^+=\frac12\Ga_2+ \frac34M$,
where $M:=\{y_4+y_3=0\}\cap S$. Then $K_S+\De^+$ is not klt at
$(0,0,1,-1)$. This implies that the singularity is not
exceptional.

\subsubsection*{}
Now we prove case by case that all singularities in
Table~\ref{table-2} are exceptional. Assuming that there exists a
regular non-klt $n$-complement $K_S+\De^+$, we derive a
contradiction or prove that $K_S+\De^+$ is klt (see
Proposition~\ref{prokhorov}). Set $\De':=\De^+-\De$. By
Corollary~\ref{support} we may assume either
$\Supp(\De^+)\not\subset\Ga$ or $\Supp(\De^+)\subset\Ga$ and
$\down{\De^+}\ne 0$. By \cite[Proposition~2]{prokhorov} (or
\cite[4.9]{Pr_lect}) and (i) of \ref{properties-Gamma}, for any
$n$-complement $K_S+\De^+$ there exists an $n$-complement
$K_{\bP}+D$ such that $D|_S=\De^+$. So we can write
\begin{equation}
\label{De}
\De^+=\sum \frac{k_i}n\Ga_i+\sum \frac{s_j}nM_j,
\end{equation}
where $k_i,s_j\in\bN\cup\{0\}$, $M_j$'s are effective (not
necessarily irreducible) curves and $M_j\sim m_jH$. We may assume
that $M_j\ne \Ga_i$, $\forall i, j$. By \ref{support} one of the
following holds $\sum \frac{s_j}nM_j\ne 0$ or $\down{\sum
\frac{k_i}n\Ga_i}\ne 0$. By the construction, $\sum
\frac{s_j}nM_j\le \De'$.

\subsubsection*{Cases $\x238{18}$, $\x2469$ and $\x239{12}$}
Then $\De'\equiv\frac1{12}H$, $\frac16H$ and $\frac16H$,
respectively. Here $\Supp(\De^+)\subset\Ga$. It is easy to check
that $\down{\De^+}=0$. This implies that $K_S+\De^+$ is klt by
Corollary~\ref{support}.

\subsubsection*{Case $\x246{10}$}
By \ref{properties-complemnets}, $\De^+\ge
\frac12\Ga_2+\frac23\Ga_3+\frac56\Ga_4$. Since $\De^+\equiv 2H$,
$2\ge \frac12+\frac23+\frac56=2$. This yields the equality
$\De^+=\frac12\Ga_2+\frac23\Ga_3+\frac56\Ga_4$. This complement is
klt by \ref{support}. Similarly, we can argue in cases
$\x334{10}$, $\x239{10}$ and $\x23{10}{12}$.

\subsubsection*{Cases $\x238{3r}$, $r\in\{3,5,7\}$,
$\x239{2r}$, $r\in\{7,8\}$} By $\De_{r=3}\le \De_{r=5}\le
\De_{r=7}$ it is sufficient to consider only case $\x2389$. Then
$\De'\equiv\frac5{12}H$. Let $M_j$ be as in \eref{De}. Then
$\frac1nM_j\le\De'$. Hence, $m_j/n\le 5/12$. This yields $n=6$ and
$m_j=2$. But then $\De^+\ge \frac56\Ga_3+\frac23\Ga_4+\frac16M$,
so $4\ge \frac56\cdot 3+\frac23\cdot 2+\frac16\cdot 2=25/6$, a
contradiction.

\subsubsection*{Case $\x3355$}
\label{3-3-5-5}
Then $\De^+\equiv H$. Thus $n\De^+$ is given by some polynomial
$p\in H^0(S,\OOO_S(n))$, where $n\in\{1,2,3,4,6\}$. On the other
hand, $H^0(S,\OOO_S(n))\ne 0$ for $n=3$ or $6$. We obtain two
cases: $\De^+=\frac13M$, where $M$ is given by $\const_3
y_3+\const_4 y_4=0$, or $\De^+=\frac16M+\frac16M'$, where $M$ and
$M'$ are given by $\const_3 y_3+\const_4 y_4=0$ and $\const_3'
y_3+\const_4' y_4=0$, respectively. Consider for example, the
first case. If $\const_3^5+\const_4^5\ne 0$, then $M$ is
irreducible. If $\const_3^5+\const_4^5=0$, then $M$ has exactly
three irreducible components. By Lemma~\ref{can-term} it is
sufficient to show that $K_S+M^{(i)}$ is plt for any irreducible
component $M^{(i)}$ of $M$. Note that $M$ contains three singular
points $\Ga_3\cap\Ga_4$ of type $\frac15(1,1)$. If $M$ is
irreducible, then $(K_S+M)\cdot M=6H^2= 2/5$ and Lemma~\ref{plt-1}
give us that $K_S+M$ is plt. If $M$ is not irreducible, then all
the irreducible components $M^{(i)}\subset M$ are smooth. By
Lemma~\ref{plt-0}, $K_S+M^{(i)}$ is plt.

\subsubsection*{Case $\x3348$}
Then $\De'\equiv\frac12H$. Obviously, $\down{\De^+}=0$. Hence,
$\Supp(\De^+)\not\subset\Ga$. We have only one possibility:
$\De^+=\frac12\Ga_4+\frac16M$, where $M:= S \cap \{y_3=\const
y_4\}$. This curve is irreducible if $1+\const^4\ne 0$ and has
exactly three irreducible components if $1+\const^4=0$. By
Lemma~\ref{can-term} it is sufficient to show that $K_S+\Ga$ is
plt for any irreducible component of $M$. Note that $M$ contains
three singular points $\Ga_3\cap\Ga_4$ of type $\frac14(1,1)$ and
$(K_S+M)\cdot M=3H^2= 1/4$. If $M$ is irreducible, then by
Lemma~\ref{plt-1}, $K_S+M$ is plt. If $M$ is not irreducible, then
all the irreducible components $M^{(i)}\subset M$ are smooth by
Lemma~\ref{plt-0}. Thus $K_S+M^{(i)}$ is plt.

\subsubsection*{Case $\x2477$}
Then $\De'\equiv\frac12H$. It is easy to see that the coefficient
of $\Ga_2$ in $\De^+$ cannot be greater than $1/2$. We have
$\De^+=\frac12\Ga_2+\frac14M$, where $M$ is cut out on $S$ by
$y_3=\const y_4=0$. If $1+\const^7\ne 0$, then $M$ is irreducible.
If $1+\const^7= 0$, then $M$ has exactly two irreducible
components. By Lemma~\ref{can-term} it is sufficient to show that
$K_S+M^{(i)}$ is plt for any irreducible component of $M$. Note
that $M$ contains two singular points $\Ga_3\cap\Ga_4$ of type
$\frac17(1,1)$. It is easy to compute also $(K_S+M)\cdot M=-4H^2=
-2/7$. If $M$ is irreducible, then by Lemma~\ref{plt-1}, $K_S+M$
is plt. If $M$ is not irreducible, then all the irreducible
components $M^{(i)}\subset M$ are smooth by Lemma~\ref{plt-0}.
Hence, $K_S+M^{(i)}$ is plt.

\subsubsection*{Case $\x3346$}
Then $\De'\equiv\frac12H$. First we claim that $K_S+\De^+$ is klt
outside of $\Ga_3$. Assume the opposite. Then $K_S+\De'$ is not
klt and $K_S+\De'+\Ga_4$ is not plt. Write
$\De'=\alpha\Ga_4+\De''$, where $\alpha\ge 0$ and $\Ga_4$ is not a
component of $\Supp(\De'')$. By Connectedness Lemma
\cite[5.7]{shokurov2}, $K_S+\De''+\Ga_4$ is not plt near $\Ga_4$.
Further, $\dif_{\Ga_4}(0)=\frac12P_1+\frac12P_2+\frac12P_3$, where
$\{P_1,P_2,P_3\}=\Ga_3\cap\Ga_4=\Sing(S)$. On the other hand,
$\De''\cdot\Ga_4\le \De'\cdot\Ga_4=1/4$. This yields
$\down{\dif_{\Ga_4}(\De'')}=0$. By Inversion of Adjunction,
$K_S+\Ga_4+\De''$ is plt near $\Ga_4$, a contradiction.

Now we claim that $K_S+\De^+$ is klt outside of $\Sing(S)$. As
above, $\dif_{\Ga_3}(0)=\frac12P_1+\frac12P_2+\frac12P_3$ and
$\De'\cdot \Ga_3=3/4$. If $K_S+\De^+$ is not klt at some point
$P\notin\Sing(S)$, then $K_S+\Ga_3+\De'$ is not plt at $P$. Thus
$\dif_{\Ga_3}(\De')\ge\frac12P_1+\frac12P_2+\frac12P_3+P$, a
contradiction.

Further, fix a point $P_1\in\Sing(S)$ (of type $A_1$), say
$P_1=(1,-1,0,0)$ and let $M$ be as in \eref{De}. We may assume
that $P_1\in M$ and $K_S+\De^+$ is not klt at $P_1$. Since
$\De'\equiv\frac12H$, there are three cases:
\begin{enumerate}
\item
$\De^+=\frac12\Ga_3+\frac14M$, where $M:=\{y_1+y_2=\const
y_4^2\}\cap S$;
\item
$\De^+=\frac12\Ga_3+\frac16\Ga_4+\frac16M$, where
$M:=\{y_1+y_2=\const y_4^2\}\cap S$;
\item
$\De^+=\frac12\Ga_3+\frac16M$, where $M:=\{y_3=y_4(\const_1
y_1+\const_2 y_2+\const_4 y_4^2)\}\cap S$.
\end{enumerate}

By Lemma~\ref{can-term} and Lemma~\ref{plt-0} it is sufficient to
show that $M$ is either smooth at $P_1$ or $M$ has two smooth
analytic components at $P_1$. Indeed, in cases (i) and (ii), $M$
is given by $y_3^2+3\const
y_4^2y_1^2-3\const^2y_4^4y_1+(\const^3+1)y_4^6=0$ in
$\bP_{y_1,y_3,y_4}(2,3,1)$. The local equation of $M$ near
$(1,0,0)$ is $x^2+3\const y^2-3\const^2y^4+(1+\const^3)y^6=0$ in
$\bC^2/\bZ_2(1,1)$. Thus, in cases (i) and (ii), $M$ has exactly
two components which are smooth curves. Similarly, in case (iii),
$M$ is smooth at $P_1$.

\subsubsection*{Case $\x2468$}
Then $\De'\equiv\frac16H$. It is easy to compute that
$\down{\De^+}=0$. Let $M_1$ be as in \eref{De}. We have
$\frac16m_1\le 1/6$. Hence, $\De'=\frac16M_1$, where
$M_1:=\{y_4=\const y_2\}\cap S$, $\const\ne 0$. If $1+\const^4\ne
0$, then $M_1\simeq\bP^1$. If $1+\const^4=0$, then $M_1$ has
exactly two components $\simeq\bP^1$. Further,
$\dif_{\Ga_4}(0)=\frac12P_1+\frac12P_2$, where
$\{P_1,P_2\}=\Ga_2\cap\Ga_4$. On the other hand,
$\De'\cdot\Ga_4=1/6$ yields $\down{\dif_{\Ga_4}(\De')}=0$. By
Inversion of Adjunction, $K_S+\Ga_4+\De'$ is plt near $\Ga_4$.
Thus $K_S+\De^+$ is klt near $\Ga_4$. Outside of $\Ga_4$ the
surface $S$ is smooth and $M_1$ has at most two irreducible
components (and they both are smooth). By Lemma~\ref{can-term},
$K_S+\De^+$ is klt.

\subsubsection*{Case $\x2466$}
Then $\De'\equiv\frac12H$. Clearly, the coefficient of $\Ga_2$ in
$\De^+$ is $1/2$. Thus $\Ga_2$ is not a component of
$\Supp(\De')$. First we claim that $K_S+\De^+$ is klt near
$\Ga_2$. Since $\De'\cdot \Ga_2=1$ and because $\Ga_2$ is
contained in the smooth locus of $S$, we have that
$\down{\dif_{\Ga_2}(\De')}$ is reduced. By Inversion of
Adjunction, $K_S+\Ga_2+\De'$ is lc near $\Ga_2$. Then
$K_S+\frac12\Ga_2+\De'$ is klt near $\Ga_2$. This proves our
claim.
\par
Since $K_S+\De^+$ is klt near $\Ga_2$, $\down{\De^+}=0$. Now we
have to show only that $K_S+\De^+$ is klt outside of $\Ga_2$.
Assume the opposite. Then $K_S+\De'$ is not klt. Write
$\De'=\alpha \Ga_3+\De''$, where $\alpha\ge 0$ and $\Ga_3$ is not
a component of $\Supp(\De')$. Then $-(K_S+\Ga_3+\De'')$ is ample.
By Connectedness Lemma \cite[5.7]{shokurov2}, $K_S+\Ga_3+\De''$ is
not plt near $\Ga_3$ and klt outside of $\Ga_3$. Note that $\Ga_3$
contains exactly two singular points $\{P_1, P_2\}=\Ga_3\cap\Ga_4$
and they are of type $\frac13(1,1)$. Thus
$\dif_{\Ga_3}(0)=\frac23P_1+\frac23P_2$. On the other hand,
$\De''\cdot \Ga_3\le \De'\cdot\Ga_3=1/3$. From this we have that
$\up{\dif_{\Ga_3}(\De'')}$ is reduced, i.e.,
$K_{\Ga_3}+\dif_{\Ga_3}(\De'')$ is lc. Again by Inversion of
Adjunction (see \cite[17.7]{uta}), $K_S+\Ga_3+\De''$ is lc near
$\Ga_3$. Therefore, $K_S+\alpha \Ga_3+\De''$ is klt near $\Ga_3$,
a contradiction.

\subsubsection*{Case $\x245{10}$}
Then $\De'\equiv\frac12H$. There are only two cases:
$\De^+=\frac12\Ga_2+\frac14M$ or
$\De^+=\frac12\Ga_2+\frac16\Ga_4+\frac16M$, where $M$ is cut out
on $S$ by $y_3=\alpha y_4^2$. It is easy to see that $M$ is
isomorphic to a (reduced) conic. So $M$ has at most two
irreducible components and they are smooth. By
Lemma~\ref{can-term}, $K_X+\De^+$ is klt outside of
$\Sing(S)=\Ga_3\cap\Ga_4$. Near $\Sing(S)$ we have $\De^+=\De'$.
It is sufficient to show that $K_S+\Ga_4+\De'$ is plt near
$\Ga_4$. Indeed, $\dif_{\Ga_4}(0)=\frac45P_1+\frac45P_2$ and
$\De'\cdot\Ga_4=1/10<1/5$. This yields
$\down{\dif_{\Ga_4}(\De')}=0$. By Inversion of Adjunction,
$K_S+\Ga_4+\De'$ is plt near $\Ga_4$.

\subsubsection*{Cases $\x2556$, $\x245{15}$, $\x2558$}
It is sufficient to consider only case $\x2556$ (when $\De$ is
smaller). Then $\De'\equiv\frac23H$. Let $\gamma$ be the
coefficient of $\Ga_4$ in $\De^+$. Assume that $\gamma>2/3$. Since
$\Ga_4\sim 5H$, we have $\gamma-2/3\le 5/6$. Taking into account
that $\gamma=k/n$, $k\in\bN$, we obtain $\gamma =3/4$ and $n=4$.
This means that $4(\De^+-\frac34\Ga_4)$, is an integral effective
divisor on $S$. On the other hand, $4\De'\equiv H$, a
contradiction. Therefore we may assume that $\gamma =2/3$ and
$n=3$ or $6$. By our assumption $\Ga_4$ is not a component of
$\De'$. First, we claim that $K_S+\Ga_4+\De'$ is plt near $\Ga_4$.
By Lemma~\ref{sing}, $S$ has on $\Ga_4$ five points
$P_1,\dots,P_5$ of type $A_1$. Thus, $\dif_{\Ga_4}(0)=\frac12\sum
P_i$. On the other hand, $\De'\cdot\Ga_4=\frac{10}3H^2=1/3$.
Therefore, $\down{\dif_{\Ga_4}(\De')}=0$. By Inversion of
Adjunction, $K_S+\Ga_4+\De'$ is plt near $\Ga_4$. Assume that
$K_S+\De^+$ is not klt outside of $\De_4$. Then $K_S+\De'$ is not
klt. Since $-(K_S+\Ga_2+\De')$ is ample, by Connectedness Lemma
\cite[5.7]{shokurov2}, $K_S+\Ga_2+\De'$ is not plt near $\Ga_2$.
Hence, $K_{\Ga_2}+\dif_{\Ga_2}(\De')$ is not klt (i.e.
$\down{\dif_{\Ga_2}(\De')}\ne 0$). On the other hand, as above,
$\dif_{\Ga_2}(0)=\frac45Q_1+\frac45Q_2$ and $\De'\cdot\Ga_2=
2/15<1/5$. So $\down{\dif_{\Ga_2}(\De')}\ne 0$, a contradiction.

\subsubsection*{Case $\x238{12}$}
Then $\De'\equiv\frac12H$. If $\Supp(\De^+)\subset\Ga$, then
$\down{\De^+}=0$ and $K_S+\De^+$ is klt by \ref{support}. If
$\Supp(\De^+)\not\subset\Ga$, then the only possibility is
$\De^+=\frac12\Ga_3+ \frac16M$, where $M=\{y_3=\const y_4^3\}\cap
S$, $\const\ne 0$. It is easy to see that
$M\simeq\{x^2+y^3+(1+\const^4)z^6=0\}\subset\bP(3,2,1)$. This
curve is irreducible and $p_a(M)=1$. If $1+\const^4\ne 0$, then
$M$ is smooth. If $1+\const^4=0$, then $M$ has a simple cusp at
$P=\{(0,0,\const,1)\}$. By Lemma~\ref{sing}, $S$ is smooth at $P$.
Further, $S$ has only singularities of type $A_1$ and
$\frac13(1,1)$ (see \ref{sing}). Lemma~\ref{plt-0} and
Lemma~\ref{can-term} give us that $K_S+\De^+$ is klt outside of
$P$. But at $P$ we have $\De^+=\frac16M$, where $M$ has a simple
cusp at $P$. At this point $K_S+\frac16M$ is also klt (see, e.g.,
\cite[8.14]{kollar}).

\clearpage \setcounter{no}{0}
\begin{table}[h]
\caption{Case: $\Pic(S)=\bZ$}
\label{table-1}
\begin{tabular}{|c|ccccp{0.7cm}p{1.9cm}p{0.1cm}p{4.5cm}|}
\hline No.&$a_1$&$a_2$&$a_3$&$a_4$&&$S$&&$\dif_S(0)$\\ \hline
\no&$2$&$3$&$7$&$r$&&$\bP^2$&&$\frac12L_1+ \frac23L_2+ \frac67L_3+
\frac{r-1}rC_1$\\
&\multicolumn{8}{|l|}{$r\in\{11,13,17,19,23,25,29,31,37,41\}$}\\
\no&$2$&$3$&$11$&{13}&&$\bP^2$&&
$\frac12L_1+\frac23L_2+\frac{10}{11}L_3+\frac{12}{13}C_1$\\
\no&$2$&$4$&$5$&$2r$&&$\bP^2$&&$\frac12L_2+\frac45C_2+\frac{r-1}rL_3$\\
&\multicolumn{6}{|l}{$r\in\{3,7,9\}$}&&\\
\no&$2$&$4$&$6$&$r$&&$\bP^2$&&$\frac12L_2+\frac23L_3+\frac{r-1}rC_2$\\
&\multicolumn{6}{|l}{$r\in\{7,11\}$}&&\\
\no&$2$&$3$&$8$&$2r$&&$\bP^2$&&$\frac23C_2+\frac34L_2+\frac{r-1}rL_3$\\
&\multicolumn{6}{|l}{$r\in\{5,7,11\}$}&&\\
\no&$2$&$3$&$10$&{14}&&$\bP^2$&&$\frac23C_2+\frac45L_2+\frac67L_3$\\
\no&$2$&$3$&$9$&{15}&&$\bP^2$&&$\frac12C_3+\frac23L_2+\frac45L_3$\\
\no&$3$&$3$&$4$&$9$&&$\bP^2$&&$\frac34C_3+\frac23L_3$\\
\no&$3$&$3$&$5$&$6$&&$\bP^2$&&$\frac45C_3+\frac12L_3$\\
\no&$3$&$4$&$4$&$4$&&$\bP^2$&&$\frac23C_4$\\
\no&$2$&$5$&$5$&$5$&&$\bP^2$&&$\frac12C_5$\\
\no&$2$&$3$&$7$&$2r$&&$\bP(1,2,1)$&&$\frac 23C_2+\frac
67L_2+\frac{r-1}rL_3$\\&\multicolumn{8}{|l|}{$r\in
\{4,5,8,10,11,13,16,17,19,20\}$}\\
\no&$2$&$3$&$8$&$r$&&$\bP(1,1,2)$&&$\frac 23C_2+\frac
34L_2+\frac{r-1}rL_3$\\&\multicolumn{6}{|l}{$r\in
\{11,13,17,19,23\}$}&&\\
\no&$2$&$3$&$8$&${20}$&&$\bP(2,1,1)$&&$\frac
23C_4+\frac12L_2+\frac 45L_3$\\
\no&$2$&$4$&$5$&$r$&&$\bP(1,1,2)$&&$\frac 12L_2+\frac
45L_3+\frac{r-1}rC_2$\\&\multicolumn{6}{|l}{$r\in
\{7,9,11,13,17,19\}$}&&\\
\no&$2$&$4$&$5$&$4r$&&$\bP(2,1,1)$&&$\frac 45C_4+\frac{r-1}rL_3$\\
&\multicolumn{6}{|l}{$r\in \{2,3,4\}$}&&\\
\no&$2$&$4$&$7$&$8$&&$\bP(2,1,1)$&&$\frac 67C_4+\frac 12L_3$\\
\no&$2$&$4$&$7$&$9$&&$\bP(1,1,2)$&&$\frac 12L_2+\frac67L_3+\frac
89C_2$\\ \no&$2$&$3$&${10}$&$r$&&$\bP(1,1,2)$&&$\frac 23C_2+\frac
45L_2+\frac{r-1}rL_3$\\ &\multicolumn{6}{|l}{$r\in \{11,13\}$}&&\\
\no&$2$&$5$&$6$&$7$&&$\bP(1,1,2)$&&$\frac 45C_2+\frac23L_2+\frac
67L_3$\\ \no&$2$&$3$&$7$&$3r$&&$\bP(3,1,1)$&&$\frac 12L_1+\frac
67C_3+\frac{r-1}rL_3$\\
&\multicolumn{6}{l}{$r\in\{3,5,9,11,13\}$}&&\\
\no&$2$&$3$&$9$&$9$&&$\bP(3,1,1)$&&$\frac12C_9$\\
\no&$2$&$3$&$9$&$r$&&$\bP(3,1,1)$&&$\frac 12L_1+\frac
23L_3+\frac{r-1}rC_3$\\&\multicolumn{6}{l}{$r\in
\{11,13,17\}$}&&\\ \no&$3$&$3$&$4$&$r$&&$\bP(1,1,3)$&&$\frac
34L_3+\frac{r-1}rC_3$\\&\multicolumn{6}{l}{$r\in \{5,7,11\}$}&&\\
\no&$3$&$3$&$5$&$7$&&$\bP(1,1,3)$&&$\frac 45L_3+\frac 67C_3$\\
\no&$2$&$5$&$6$&$6$&&$\bP(3,1,1)$&&$\frac 45C_6$\\ \hline
\end{tabular}
\end{table}

\begin{table}[h]
\begin{tabular}{|c|ccccp{0.6cm}p{1.9cm}p{0.1cm}p{4.5cm}|}
%\begin{tabular}{|c|ccccp{0.1cm}p{2cm}p{0.1cm}p{4.7cm}|}
\hline
No.&$a_1$&$a_2$&$a_3$&$a_4$&&\multicolumn{1}{c}{$S$}&&$\dif_S(0)$\\
\hline \no&$2$&$3$&$8$&$8r$&&$\bP(4,1,1)$&&$\frac
23C_8+\frac{r-1}rL_3$\\ &\multicolumn{6}{l}{$r\in \{1,2\}$}&&\\
\no&$3$&$4$&$4$&$5$&&$\bP(4,1,1)$&&$\frac 23L_1+\frac 45C_4$\\
\no&$2$&$5$&$5$&$r$&&$\bP(5,1,1)$&&$\frac 12L_1+\frac{r-1}rC_5$\\
&\multicolumn{6}{l}{$r\in \{7,9\}$}&&\\
\no&$2$&$3$&$10$&$10$&&$\bP(5,1,1)$&&$\frac 23C_{10}$\\
\no&$2$&$3$&$7$&{35}&&$\bP(7,1,1)$&&$\frac 12C_7+\frac 23L_1+\frac
45L_3$\\ \no&$2$&$3$&$11$&$11$&&$\bP(11,1,1)$&&$\frac
12C_{11}+\frac 23L_1$\\
\no&$2$&$3$&$7$&${6r}$&&$\bP(3,2,1)$&&$\frac
67C_6+\frac{r-1}rL_3$\\
&\multicolumn{6}{|l}{$r\in\{2,3,4,5,6\}$}&&\\
\no&$2$&$3$&${11}$&${12}$&&$\bP(3,2,1)$&&$\frac{10}{11}C_6+\frac
12L_3$\\ \no&$2$&$3$&$7$&${14r}$&&$\bP(7,2,1)$&&$\frac
23C_{14}+\frac{r-1}rL_3$\\&\multicolumn{6}{|l}{$r\in\{1,2\}$}&&\\
\no&$2$&$3$&$7$&${21}$&&$\bP(7,3,1)$&&$\frac 12C_{21}$\\ \hline
\end{tabular}
\end{table}

\begin{table}[ht]
\caption{Case: $\rho(S)>1$}
\label{table-2}
\begin{tabular}{|c|ccccccp{2.4cm}|}
\hline No.&${a_1}$&${a_2}$&${a_3}$&${a_4}$&\multicolumn{2}{c}
{$\bP(\bar{\bp})\supset S$}&$\dif_S(0)$\\ \hline
\no&$2$&$5$&$5$&$2r$&$\bP(5,2,2,5)$&$\y2552$&$\frac{r-1}r\Ga_4$\\
&\multicolumn{5}{l}{$r\in \{3,4\}$}&&\\
\no&$2$&$4$&$5$&${15}$&$\bP(5,5,2,2)$&$\y2255$&$\frac12\Ga_2+\frac23\Ga_4$\\
\no&$2$&$3$&$8$&$3r$&$\bP(3,2,3,2)$&$\y2323$&$\frac34\Ga_3+\frac{r-1}r\Ga_4$\\
&\multicolumn{5}{l}{$r\in \{3,5,7\}$}&&\\
\no&$2$&$3$&$9$&$2r$&$\bP(3,2,2,3)$&$\y2332$&$\frac23\Ga_3+\frac{r-1}r\Ga_4$\\
&\multicolumn{5}{l}{$r\in \{5,7,8\}$}&&\\
\no&$3$&$3$&$4$&${10}$&$\bP(2,2,3,3)$&$\y3322$&$\frac12\Ga_3+\frac45\Ga_4$\\
\no&$2$&$4$&$5$&${10}$&$\bP(5,5,2,1)$&$\y225{10}$&$\frac12\Ga_2$\\
\no&$2$&$3$&$8$&${18}$&$\bP(3,2,3,1)$&$\y2326$&$\frac34\Ga_3+\frac23\Ga_4$\\
\no&$2$&$3$&${10}$&${12}$&$\bP(3,2,3,1)$&$\y2326$&$\frac45\Ga_3+\frac12\Ga_4$\\
\no&$2$&$4$&$6$&$9$&$\bP(3,3,1,2)$&$\y2263$&$\frac12\Ga_2+\frac23\Ga_4$\\
\no&$2$&$4$&$6$&${10}$&$\bP^3$&$\y2222$&$\frac12\Ga_2+\frac23\Ga_3+\frac45\Ga_4$\\
\no&$2$&$3$&$8$&${12}$&$\bP(6,4,3,1)$&$\y234{12}$&$\frac12\Ga_3$\\
\no&$3$&$3$&$4$&$6$&$\bP(2,2,3,1)$&$\y3326$&$\frac12\Ga_3$\\
\no&$2$&$3$&$9$&${12}$&$\bP(3,2,2,1)$&$\y2336$&$\frac23\Ga_3+\frac12\Ga_4$\\
\no&$2$&$4$&$7$&$7$&$\bP(7,7,2,2)$&$\y2277$&$\frac12\Ga_2$\\
\no&$3$&$3$&$4$&$8$&$\bP(4,4,3,3)$&$\y3344$&$\frac12\Ga_4$\\
\no&$3$&$3$&$5$&$5$&$\bP(5,5,3,3)$&$\y3355$&$0$\\
\no&$2$&$4$&$6$&$8$&$\bP(2,1,2,1)$&$\y2424$&$\frac23\Ga_3+\frac12\Ga_4$\\
\no&$2$&$4$&$6$&$6$&$\bP(3,3,1,1)$&$\y2266$&$\frac12\Ga_2$\\
\hline
\end{tabular}
\end{table}
\clearpage
\subsubsection{}
\label{Q-factorial}
One can check that all surfaces $S$ in Table~\ref{table-2} have
the Picard number $\rho(S)>1$. Hence, all these singularities are
not analytically $\bQ$-factorial:

\begin{proposition*}
\label{non-Q-fact}
Let $(X,P)$ be an analytic germ of a klt singularity and let
$f\colon (Y,S)\to X$ be a plt blow-up such that $f(S)=P$. Then
\begin{enumerate}
\item
$\Pic(Y)=H^2(Y,\bZ)=H^2(S,\bZ)=\Pic(S)$;
\item
if $(X,P)$ is analytically $\bQ$-factorial, then $\Pic(S)=\bZ$.
\end{enumerate}
\end{proposition*}
\begin{proof}
We have an exact sequence
\[
0\longrightarrow\bZ\longrightarrow \OOO_Y
\stackrel{\exp}{\longrightarrow} \OOO_Y^*\longrightarrow 0.
\]
By Kawamata-Viehweg vanishing $R^if^*\OOO_Y=0$, $i>0$. Hence,
$\Pic(Y)=H^2(Y,\bZ)$. Similarly, $H^2(S,\bZ)=\Pic(S)$. Since $Y$
is an analytic germ near $S$, $H^2(Y,\bZ)=H^2(S,\bZ)$. If $(X,P)$
is analytically $\bQ$-factorial, then $\rho(Y/X)=1$ and
$\rk\Pic(S)=1$.
\end{proof}

\begin{remark*}
(i) In case $\x246{10}$, $S\simeq\bP^1\times\bP^1$. In cases
$\x3346$ and $\x239{12}$, the projection $S\to\bP(2,2,1)=\bP^2$ is
a double cover ramified along $\{(x^3+y^3+z^3)z=0\}\subset\bP^2$.
Hence, $S$ is a Gorenstein del Pezzo surface of degree $2$ having
exactly three singular points which are of type $A_1$. In case
$\x2468$, $\bP(2,1,2,1)$ is isomorphic to a cone over a conic in
$\bP^4$. Here $S$ is isomorphic to an intersection of two quadrics
in $\bP^4$. In cases $\x238{18}$, $\x23{10}{12}$ and $\x2469$, the
surface $S$ has exactly two singular points of type $A_2$. It is a
Gorenstein del Pezzo surface of degree $3$, a cubic in $\bP^3$.
\par
(ii) In cases \x{2}{3}{7}{r},
$r\in\{11,13,17,19,23,25,29,31,37,41\}$ and \x{2}{3}{11}{13}, the
pairs $(S,\De)$ have Shokurov's invariant $\delta(S,\De)=2$ (see
\cite[\S 5]{shokurov}, and also \cite{Pr_lect}). The singularity
$\x2478$ and all singularities with $S\simeq\bP(1,2,3)$ give us a
log-del Pezzo surface $(S,\De)$ of the so-called ``elliptic type''
(see \cite[Table~1, No.~1]{Abe}).
\end{remark*}
%%%%%%%%%%%%%%%%%%%%%%%%%%%%%%%%%%%%%%%%%%%%%%%%%%%%%%%

\end{document}